\documentclass[11pt]{article}
\usepackage[T1]{fontenc}
\usepackage{lmodern}
\usepackage[margin=1.05in]{geometry}
\usepackage{amsmath,amssymb,amsthm,mathtools}
\usepackage{booktabs,array,microtype}
\usepackage[hidelinks]{hyperref}
\hypersetup{pdftitle={Fixed-point-free involution factorizations (and the permutations that admit them).},pdfauthor={Charles Burnette},pdfsubject={Fixed-point-free involution factorizations, discrete limits, medians, and moments}}
\usepackage{enumitem}
\setlist[enumerate]{itemsep=3pt,topsep=5pt}
\newcommand{\invol}{\mathsf{invol}}
\newcommand{\E}{\mathbb E}
\renewcommand{\P}{\mathbb P}
\newcommand{\Var}{\operatorname{Var}}

\newcommand{\Li}{\operatorname{Li}}

\newcommand{\dd}{\,\mathrm d}
\newcommand{\ee}{\mathrm e}
\newcommand{\TV}{\mathrm{TV}}
\newtheorem{theorem}{Theorem}[section]
\newtheorem{lemma}[theorem]{Lemma}
\newtheorem{proposition}[theorem]{Proposition}
\newtheorem{corollary}[theorem]{Corollary}
\theoremstyle{definition}
\newtheorem{definition}[theorem]{Definition}
\newtheorem{example}[theorem]{Example}
\theoremstyle{remark}
\newtheorem{remark}[theorem]{Remark}
\numberwithin{equation}{section}
\allowdisplaybreaks[1]
\title{Fixed-point-free involution factorizations\\
(and the permutations that admit them)}
\author{Charles Burnette\\[3pt]
\small Department of Mathematics\\
\small Xavier University of Louisiana\\
\small New Orleans, LA 70125-1098\\
\small \texttt{cburnet2@xula.edu}}
\date{}
\begin{document}
\maketitle
\begin{abstract}
A permutation is a product of two fixed-point-free involutions if and only if all its cycle multiplicities are even. We call such permutations \emph{acme} and study their number of fixed-point-free involution factorizations. An explicit coefficient expansion gives the proportion of the class of acme permutations within the symmetric group. For an acme permutation chosen uniformly at random from the permutations of degree $2n$, we remove the two cycles whose common length exceeds $n/2$, when such cycles are present. The numbers of cycles of each length that remain, divided by two, converge jointly to independent random variables. We obtain a sharp total variation expansion for the convergence rate. This yields an explicit discrete limit law for the factorization count divided by $n$. The factorization count has unique median $n$ for every positive $n$ except $n=3$, whose median is $2$. Its moments display several distinct regimes: fractional moments below order 1 converge on the typical scale, the first moment has order $n^{3/2}$, each fixed order between 1 and 2 admits a full saddle-point equivalent, the second contains a partition-type exponential factor, and every fixed moment above order 2 is asymptotically carried by the identity permutation. The proofs combine cycle pairing, the Shepp--Lloyd model, hybrid Darboux--singularity analysis, and saddle-point estimates with explicit control of the complementary arcs.
\end{abstract}
\noindent\textit{2020 Mathematics Subject Classification.} 05A05, 05A16, 60C05, 60F05.
\par\smallskip
\noindent\textit{Keywords.} Involution factorization; random permutation; conditioned cycle counts; singularity analysis; median; moments; long cycles.

\section{Introduction}
\subsection{Background}
An involution is a permutation that is equal to its own inverse. If $\sigma$ is a permutation in the symmetric group $\mathfrak{S}_d$, the statistic
\[
 \invol(\sigma):=\#\{(\tau_1,\tau_2)\in \mathfrak{S}_d\times \mathfrak{S}_d:
 \tau_1^2=\tau_2^2=\text{id},\ \sigma=\tau_2\tau_1\}
\]
counts the ordered decompositions of $\sigma$ as a product of two involutions. Contemporary work on the combinatorics of involution factorizations originates in~\cite{Lugo} and~\cite{PetersenTenner}. Moreover, Lugo conjectured in~\cite{Lugo} that $\invol$ is asymptotically lognormal for permutations drawn uniformly at random. This was subsequently proven in~\cite{BurnetteSchmutz} and strengthened into a functional central limit theorem for $\log(\invol)$ over Ewens Sampling Formula permutations generated by the Feller coupling~\cite{BurnetteEwens}.

The present paper studies the fixed-point-free analogue of $\invol$, which we denote by $\invol^*$. Specifically, if $\mathcal{D}_d$ denotes the set of all derangements in $\mathfrak{S}_d$, then
\[
 \invol^*(\sigma):=\#\{(\tau_1,\tau_2)\in \mathcal{D}_d\times \mathcal{D}_d:
 \tau_1^2=\tau_2^2=\text{id},\ \sigma=\tau_2\tau_1\}.
\]
There are several reasons to care about this specialized statistic. Roberts and Vivaldi~\cite{RobertsVivaldi2005,RobertsVivaldi2009} proposed products of involutions as a means to explicate time-reversal symmetry in dynamical systems with finite phase space. They emphasized that the number of fixed points in the involutive factors is key to categorizing the dynamics of a reversible map. A study of fixed-point-free involution factorizations is therefore a natural first step towards understanding the boundary case in which both fixed sets are empty.

Tout, Goupil, and Al-Ayyoub~\cite{ToutGoupilAlAyyoub25} have also connected this statistic with the structure coefficients of the center of the symmetric group algebra. For a partition $\lambda$ of $d$, let $C_\lambda$ be the conjugacy class of that cycle type and let $K_\lambda=\sum_{\rho\in C_\lambda}\rho$ be its class sum. The identity
\[
 K_\lambda K_\mu=\sum_{\nu\vdash d}c^\nu_{\lambda,\mu}K_\nu
\]
defines the connection coefficients. If $\sigma\in \mathfrak{S}_{2n}$ has cycle type $\nu$, then
$\invol^*(\sigma)=c^\nu_{2^n,2^n}$: the coefficient counts factorizations of one fixed permutation of type $\nu$.

Earlier work of Moran~\cite{Moran76a,Moran76b} treats bireflections and products of prescribed involution classes; related questions concern several conjugate factors~\cite{Moran75,Moran76c}. Mikhailov~\cite{Mikhailov17} studies logarithmic limit laws with a specified cycle in one factor, while Malyshev~\cite{Malyshev24a,Malyshev24b} considers products of three and four fixed-point-free involutions. Applications include permutation algorithms based on simultaneous swaps~\cite{YangEllisMamakaniRuskey13} and algebraic models of genome rearrangement~\cite{BhatiaEgriNagyFrancis15}. Our chief concern here is the distribution of $\mathsf{invol}^*$ over uniform random ``admissible'' permutations.

\subsection{The admissible permutations}\label{sec:exact}
Since fixed-point-free involutions are products of disjoint transpositions, it is useful to interpret them as perfect matchings of $[2n]=\{1,2,\ldots,2n\}$. When one overlays two perfect matchings $\tau_1$ and $\tau_2$ of $[2n]$, say with one matching colored red and the other blue, the resulting bicolored multigraph is $2$-regular, with every vertex incident to an edge of each color. Consequently, $\tau_1\cup\tau_2$ is a disjoint union of even cycles with alternately colored edges. A common matching edge is retained twice and forms a component on two vertices. Reading the product around a component of length $2k$ produces two disjoint cycles of length $k$; see~\cite{Lugo,ToutGoupilAlAyyoub25}. Therefore every cycle length appears an even number of times in $\tau_2\tau_1$.

\begin{definition}\label{def:acme}
A permutation $\sigma\in \mathfrak{S}_{2n}$ is called \emph{acme}, short for ``\textbf{a}ll \textbf{c}ycle \textbf{m}ultiplicities \textbf{e}ven,'' if $c_k$ is even for every $k\geq1$, where $c_k := c_k(\sigma)$ denotes the number of its $k$-cycles that $\sigma$ has. We write $\mathcal E_{2n}$ for the set of these permutations and $a_n:=|\mathcal E_{2n}|/(2n)!$ for its proportion in $\mathfrak{S}_{2n}$.
\end{definition}

As outlined above, being an acme permutation is a necessary condition to admit a fixed-point-free involution factorization. Sufficiency, together with an explicit expression for $\invol^*$, is established through the following results collected from the research literature.

\begin{definition}\label{def:reflections}
For each integer $d\geq1$ and $a\in\mathbb Z/d\mathbb Z$, define $I_a^{(d)} : \mathbb Z/d\mathbb Z \rightarrow \mathbb{Z}/d\mathbb{Z}$ to be the function given by $I_a^{(d)}(x)=a-x\pmod d,$ and define $J_a^{(d)}:\mathbb{Z}/2d\mathbb{Z} \rightarrow \mathbb{Z}/2d\mathbb{Z}$ to be the involution whose $d$ transpositions are $(x, d + I_a^{(d)}(x)),$ $x = 0, 1, \ldots, d - 1.$
\end{definition}

\begin{lemma}[Lugo~\cite{LugoThesis}; Yang--Ellis--Mamakani--Ruskey~\cite{YangEllisMamakaniRuskey13}; Burnette--Schmutz~\cite{BurnetteSchmutz}]\label{lem:single-cycle}
The $d$-cycle $x\mapsto x+1$ on $\mathbb Z/d\mathbb Z$ has exactly $d$ involution factorizations. They are $I_a^{(d)} \circ I_{a-1}^{(d)}, a\in\mathbb Z/d\mathbb Z.$ In each such factorization at least one factor has a fixed point.
\end{lemma}

Lemma~\ref{lem:single-cycle} is a specific instance of the Goupil--Schaeffer formula for decompositions of a full cycle into two permutations of prescribed cycle types~\cite{GoupilSchaeffer98}; see also Lugo~\cite[Lemma~5.7.2]{LugoThesis}. We record the formula to place the reflection construction in its class-algebra setting.\label{sec:class-algebra}

\begin{proposition}[Goupil--Schaeffer \cite{GoupilSchaeffer98}]
\label{prop:goupil-schaeffer}
Let $\lambda=(\lambda_1,\ldots,\lambda_\ell)$ and
$\mu=(\mu_1,\ldots,\mu_r)$ be partitions of $d$. If
$g=(d+1-\ell-r)/2$ is not a nonnegative integer, then
$c^{(d)}_{\lambda,\mu}=0$. Otherwise, write
$z_\lambda=\prod_{j\geq1}j^{\alpha_j}\alpha_j!$, where $\alpha_j$ is
the multiplicity of the part $j$ in $\lambda$, and define $z_\mu$ similarly.
For every nonnegative integer $v$, put
\[
 U_v(\lambda)=
 \sum_{\substack{i_1,\ldots,i_\ell\geq0\\i_1+\cdots+i_\ell=v}}
 \prod_{a=1}^\ell\binom{\lambda_a}{2i_a+1}.
\]
Then
\[
 c^{(d)}_{\lambda,\mu}
 =\frac{d}{2^{2g}z_\lambda z_\mu}
 \sum_{g_1+g_2=g}
 (\ell+2g_1-1)!(r+2g_2-1)!
 U_{g_1}(\lambda)U_{g_2}(\mu),
\]
where $g_1,g_2$.
\end{proposition}

For involution cycle types, every part is $1$ or $2$, so $U_v$ vanishes
unless $v=0$. If $d=2q+1$, each factor of a $d$-cycle has type
$1\,2^q$, and the formula gives $d$ factorizations. If $d=2q$, the two
types are $2^q$ and $1^2 2^{q-1}$, in either order, and the formula gives
$q$ factorizations in each order. These are exactly the reflection
factorizations above. The general formula provides useful context, while
the elementary orbit argument supplies everything required for the
fixed-point-free product formula.

\begin{lemma}[Burnette--Schmutz~\cite{BurnetteSchmutz}]\label{lem:cycle-orbits}
Suppose $\sigma=\tau_2\tau_1$ is an involution factorization and $\mathcal{O}$ is a cycle orbit of $\sigma$. Then $\tau_1(\mathcal{O})=\tau_2(\mathcal{O})$ is a cycle orbit of $\sigma$ of length $|\mathcal{O}|$.
\end{lemma}

\begin{definition}\label{def:exchange}
Let $\mathcal{O}_1$ and $\mathcal{O}_2$ be distinct cycle orbits of $\sigma$. An involution factorization $\sigma=\tau_2\tau_1$ \emph{exchanges} $\mathcal{O}_1$ and $\mathcal{O}_2$ if $\tau_1(\mathcal{O}_1)=\tau_2(\mathcal{O}_1)=\mathcal{O}_2$.
\end{definition}
Thus every cycle of a product of two involutions is either invariant under both factors or exchanged with another cycle of the same length.

\begin{lemma}[Lugo~\cite{LugoThesis}; Burnette--Schmutz~\cite{BurnetteSchmutz}]\label{lem:exchange-cycles}
If $\sigma = (0, 1, \ldots, d - 1)(d, d+1, \ldots dn-1),$ then there are precisely $d$ ways to write $\sigma$ as a product of two involutions of $\{0, 1, \ldots, 2d-1\}$ that exchange the two cycles of $\sigma.$ The $d$ factorizations are of the form $\sigma = J_a^{(d)} \circ J_{a - 1}^{(d)},$ $a \in \mathbb{Z}/d\mathbb{Z}.$
\end{lemma}

We now derive the formula for $\invol^*$ by a self-contained argument using the canonical action of the factors on the cycles of their product.

\begin{theorem}[Lugo~\cite{LugoThesis}; Tout, Goupil, and Al-Ayyoub~\cite{ToutGoupilAlAyyoub25}]\label{thm:formula}
A permutation $\sigma\in\mathfrak{S}_{2n}$ is a product of two fixed-point-free involutions if and only if every $\sigma$ is acme. In that case,
\begin{equation}\label{eq:formula}
\invol^*(\sigma)=\prod_{k\ge1}(c_k(\sigma)-1)!!\,k^{c_k(\sigma)/2},
\end{equation}
where we adopt the convention $(-1)!!=1$.
\end{theorem}
\begin{proof}
By Lemma~\ref{lem:cycle-orbits}, every cycle is either invariant or exchanged with another cycle of the same length under any involution factorization. Lemma~\ref{lem:single-cycle} rules out an invariant cycle when both factors are fixed-point-free. Conversely, Lemma~\ref{lem:exchange-cycles} shows that exchanging every cycle with a distinct cycle of the same length produces two fixed-point-free involutions. Such a pairing is possible exactly when every $c_k$ is even.

There are $(c_k-1)!!$ ways to pair the $k$-cycles, and each pair has $k$ possible factorizations. This contributes a factor of $(c_k-1)!!\,k^{c_k/2}$ to $\mathsf{invol}^*(\sigma).$For $c_k=0$ the contribution is $1$. Multiplication across all $k$ gives the formula.
\end{proof}

The next section uses the hybrid method of Flajolet, Fusy, Gourdon, Panario, and Pouyanne~\cite{FFGPP} to asymptotically enumerate the acme permutations. We then study the factorization count of a uniform acme permutation. The principal distinction from the unrestricted problem is that the limiting fluctuations are discrete: two long cycles supply the main factor, while the numbers of cycles of each remaining length supply the random part. This is an instance of the convergent regime of conditioned combinatorial structures~\cite{BarbourGranovsky}, for which we give a direct proof. The later sections determine the sharp median, the moment asymptotics in every fixed-order regime, and the rate of convergence of the cycle counts left after removing the two long cycles.

\section{Asymptotic enumeration of acme permutations}\label{sec:enumeration}
\subsection{The exponential generating function}
Choose $\sigma_{2n}$ uniformly from $\mathcal E_{2n}$, meaning that each acme permutation has the same probability, and write $M_{k,n}=c_k(\sigma_{2n})/2$. Then
\begin{equation}\label{eq:profilelaw}
w(\omega)=\prod_{k\ge1}\frac1{k^{2m_k}(2m_k)!},\qquad
a_n=\sum_{|\omega|=n}w(\omega),\qquad
\P((M_{k,n})=\omega)=\frac{w(\omega)}{a_n}
\end{equation}
when $|\omega|=n$.
Due to Cauchy's Formula for the number of permutations with a given cycle type, the generating function for acme permutations, counted according to the size of the ground set (indicated by the variable $z$) and the number of $k$-cycles (indicated by the variable $x_k$), exponential in $z$ and ordinary in the $x_k,$ is
\begin{align}
\mathbf{E}(\mathbf{x}, z) &:= \sum_{\substack{\lambda \vdash n, \\ \text{all}\ \lambda_k\ \text{even}}} z^n\prod_{k = 1}^{n} \left(\frac{x_k}{k}\right)^{\lambda_k}\frac{1}{\lambda_k!} = \prod_{k = 1}^{\infty} \sum_{\text{even}\ \lambda_k} \left(\frac{x_kz^k}{k}\right)^{\lambda_k}\frac{1}{\lambda_k!}
\\ &= \prod_{k = 1}^{\infty} \cosh\!\left(\frac{x_kz^k}{k}\right).\label{eq:multivar}
\end{align}
Its univariate specialization is the generating function for these permutations wholesale:
\begin{align}
E(z)&:=\sum_{n=0}^{\infty}|\mathcal E_{2n}|\frac{z^{2n}}{(2n)!}=E(1,1,\ldots;z)=F(z^2),\notag\\
F(z)&:=\prod_{k=1}^{\infty}\sum_{r=0}^{\infty}\frac{z^{kr}}{k^{2r}(2r)!}.\label{eq:generating}
\end{align}
Thus $[z^{2n}]E(z)=a_n$, while the odd coefficients of $E(z)$ vanish. The constant term is $1$, corresponding to the empty permutation. The power series starts off as
\begin{equation}\label{eq:initial-counts}
E(z)=1+\frac{z^2}{2!}+4\frac{z^4}{4!}+86\frac{z^6}{6!}
        +2\,696\frac{z^8}{8!}
       +168\,232\frac{z^{10}}{10!}
        +15\,948\,032\frac{z^{12}}{12!}+\cdots.
\end{equation}
These counts are listed as sequence A130268 in the OEIS~\cite{OEIS}. The coefficient $[z^d]E(z)$ is the probability that a uniform permutation of $[d]$ is acme. The leading asymptotic is also recorded in that entry. We now establish it, together with its first logarithmic correction. First, we determine the radius of convergence.

\begin{lemma}\label{lem:radius}
Both $E(z)$ and $F(z)$ have radius of convergence $1$.
\end{lemma}
\begin{proof}
Let $B_{2n}$ be the number of permutations in $\mathfrak{S}_{2n}$ consisting of two disjoint $n$-cycles. Since such permutations are acme, it follows that
    \begin{equation}
        \frac{B_{2n}}{(2n)!} \leq \frac{a_{2n}}{(2n)!} \leq 1.
    \end{equation}
    Now note that
    \begin{equation}
        B_{2n} = \frac{1}{2}\binom{2n}{n}(n-1)!^2 = \frac{(2n)!}{2n^2}.
    \end{equation}
    Hence
    \begin{equation}
        \lim_{n \to \infty} \sqrt[2n]{\frac{B_{2n}}{(2n)!}} = \lim_{n \to \infty} \sqrt[2n]{\frac{1}{2n^2}} = 1. 
    \end{equation}
    The result now follows from the Squeeze Theorem and the Cauchy-Hadamard Formula.
\end{proof}

\subsection{The exp-log schema}
We first introduce the normalized tangent numbers $t_j$, the coefficients in
\[
 \tan u=\sum_{j=0}^{\infty} t_j u^{2j+1},\qquad |u|<\pi/2.
\]
Integration gives
\[
 -\log\cos u=\sum_{j=1}^{\infty}t_{j-1}\frac{u^{2j}}{2j},
\]
where we have taken the principal branch of the logarithm. Hence
\begin{equation}\label{eq:bj-definition}
 \log\cosh u=\log\cos(iu)=\sum_{j=1}^{\infty}b_ju^{2j},\qquad
 b_j=\frac{(-1)^{j-1}t_{j-1}}{2j}.
\end{equation}
Writing $\Li_s(z)=\sum_{k=1}^{\infty} z^k/k^s$ in the unit disc, the exponential-logarithmic reorganization of $E$ gives
\begin{align}
 E(z)
 &=\exp\left(\sum_{k=1}^{\infty}\log\cosh(z^k/k)\right)\notag\\
 &=\exp\left(\sum_{k=1}^{\infty}\sum_{j=1}^{\infty}b_j\frac{z^{2kj}}{k^{2j}}\right)
 =\exp\left(\sum_{j=1}^{\infty}b_j\Li_{2j}(z^{2j})\right).
 \label{eq:E-explog}
\end{align}
All summation rearrangements are justified by absolute convergence, even on the closed unit disc. Indeed, $\log\cosh u$ is analytic for $|u|<\pi/2$, so for every fixed $1<\rho<\pi/2$ there is a constant $B_\rho$ such that
\begin{equation}\label{eq:bj-bound}
 |b_j|\leq B_\rho\rho^{-2j}.
\end{equation}
Moreover, recall that $\text{Li}_{\nu}(1) = \zeta(\nu)$ monotonically decreases to 1 in the limit as $\nu \to \infty.$ Thus, by Dirichlet's condition, $\sum_{j\geq1}b_j[\zeta(2j)-1]$ and, subsequently, $\sum_{j\geq1}b_j\zeta(2j)$ both converge absolutely. An appeal to Abel's theorem then provides a fast convergent series for $E(1).$ Put
\begin{equation}\label{eq:G-definition}
 G:=\sum_{j=1}^{\infty}b_j\zeta(2j),\qquad
 C:=e^G=E(1)=F(1)=\prod_{k=1}^{\infty}\cosh(1/k)
   =2.1164655365054847758\ldots.
\end{equation}
In particular, we can write
\begin{equation}\label{eq:E-centered}
 E(z)=\ee^G\exp\left(\sum_{j=1}^{\infty}b_j
       [\Li_{2j}(z^{2j})-\zeta(2j)]\right).
\end{equation}
The half-degree series
\begin{equation}\label{eq:F-definition}
 F(z)=\sum_{n=0}^{\infty}a_nz^n
 =\prod_{k=1}^{\infty}\sum_{m=0}^{\infty}\frac{z^{km}}{k^{2m}(2m)!}
\end{equation}
has the corresponding expansion
\begin{equation}\label{eq:F-polylogs}
 \log F(z)=\frac12\Li_2(z)+R(z),\qquad
 R(z)=\sum_{j=2}^{\infty}b_j\Li_{2j}(z^j).
\end{equation}
For the correction term, define
\begin{equation}\label{eq:kappa-definition}
 h=R'(1)=\sum_{j=2}^{\infty}jb_j\zeta(2j-1),\qquad
 \kappa=\gamma-\frac12+2h=-0.2240294326559579230\ldots,
\end{equation}
where $\gamma$ is Euler's constant. The geometric bound above gives absolute convergence of this series as well.

\subsection{Proof of the asymptotic estimate}
In this section, we will derive the leading term asymptotic of $F(z),$ including the smoothness needed to control the coefficient error. The second term also illustrates why an analytic term in a local expansion, although irrelevant to the first approximation, cannot be discarded before multiplying expansions.

\begin{theorem}\label{thm:enumeration}
As $n\to\infty$,
\begin{equation}\label{eq:an-two-terms}
 a_n=\frac{C}{2n^2}
 \left(1+\frac{\log n+\kappa}{n}
       +O\left(\frac{(\log n)^2}{n^2}\right)\right).
\end{equation}
In particular, a uniform random permutation of $[2n]$ is acme with probability asymptotic to $C/(2n^2)$. In odd degree this probability is zero.
\end{theorem}

\subsubsection{Darboux's method and singularity analysis}

The proof uses the hybrid method of Flajolet, Fusy, Gourdon, Panario, and Pouyanne \cite{FFGPP}. For completeness, we give a bounded version sufficient for the present product. Appendix~\ref{app:singularity} collects the underlying definitions and coefficient formulas, and explains how the two dominant singularities appear in the full-degree series $E(z)$. This avoids assuming that the whole product continues analytically through the unit circle.

A function analytic in the open unit disc $\mathbb{D}$ is called $C^s$ on $\mathbb{D}$ if its derivatives of orders $0,\ldots,s$ extend continuously to $\overline{\mathbb{D}}.$

\begin{lemma}[Darboux's estimate]\label{lem:darboux}
If $H$ is analytic in $\mathbb{D}$ and $C^s$ on $\overline{\mathbb{D}}$, where $s$ is a nonnegative integer, then $[z^n]H(z)=o(n^{-s}).$
\end{lemma}
\begin{proof}
The boundary function $\theta\mapsto H(e^{i\theta})$ is $C^s$ and periodic. Cauchy's formula, followed by $s$ integrations by parts, identifies its $n^{\text{th}}$ Fourier coefficient as $(in)^{-s}$ times the $n^{\text{th}}$ Fourier coefficient of its $s^{\text{th}}$ angular derivative. The latter tends to zero by the Riemann--Lebesgue lemma. For $s=0$ the same conclusion follows without integration by parts.
\end{proof}

By a log-power expansion at a point $\xi$ of the unit circle we mean a finite expansion in terms
\[
 (1-z/\xi)^\alpha\bigl(\log(1-z/\xi)\bigr)^r,
 \qquad r\in\{0,1,2,\ldots\},
\]
with specified remainder, using branches in a slit neighborhood of $\xi$. An indented neighborhood of the closed disk excludes a small outward-pointing wedge at each of finitely many such points. The angle of each excluded wedge is fixed and less than $\pi$.

\begin{lemma}[Bounded hybrid transfer]\label{lem:hybrid}
Suppose $F=PQ$, where $P,Q$ are analytic in the open unit disk. Let $Z$ be a finite subset of the unit circle and let $m\geq1$ be an integer. Assume that:
\begin{enumerate}
\item $Q$ is $C^{2m}$ on the closed disk;
\item $P$ continues to an indented neighborhood of the closed disk with exceptional points in $Z$, and is bounded in a sufficiently small indented neighborhood of every point of $Z$;
\item at every point of $Z$, the function $P$ has log-power expansions to every prescribed finite order.
\end{enumerate}
Then $[z^n]F(z)$ is obtained, with an error $o(n^{-m})$, by transferring the finitely many local log-power terms in the product of the expansion of $P$ and the Taylor polynomial of $Q$ through degree $m$ at each point of $Z$. One may retain additional local terms whose transferred coefficients are $o(n^{-m})$.
\end{lemma}
\begin{proof}
Choose a polynomial $H$ whose derivatives of orders $0,\ldots,m$ agree with those of $Q$ at every point of $Z$. Such a polynomial exists by Hermite interpolation. Since $Q$ has $m+1$ continuous derivatives on the closed disk, for $0\leq r\leq m$ and $z\to\xi$ inside the disk,
\[
 (Q-H)^{(r)}(z)=O\bigl(|z-\xi|^{m+1-r}\bigr).
\]
The bound follows directly from Taylor's formula along the segment joining $\xi$ to $z$. The indented continuation and local boundedness of $P$, together with Cauchy's derivative estimates on circles of radius proportional to $|z-\xi|$, give
\[
 P^{(r)}(z)=O\bigl(|z-\xi|^{-r}\bigr).
\]
Leibniz's formula therefore shows that every derivative of $P(Q-H)$ of order at most $m$ tends to zero at each $\xi\in Z$. Away from $Z$ the required continuous extensions are immediate. Hence $P(Q-H)$ is $C^m$ on the closed disk, and Lemma~\ref{lem:darboux} gives
\[
 [z^n]F(z)=[z^n](P(z)H(z))+o(n^{-m}).
\]
The function $PH$ has an indented continuation with only finitely many singularities. The multiple-singularity form of the transfer theorem \cite{FO,FFGPP} applies after its local expansions have been taken to sufficiently high order. At each $\xi$, the polynomial $H$ agrees with the Taylor polynomial of $Q$ through degree $m$. Their difference, multiplied by the bounded function $P$, has order at least $m+1$ there and contributes only terms below the required coefficient precision. This proves the assertion.
\end{proof}

The coefficient calculations themselves can be made exact. For fixed $\alpha$ and sufficiently large $n$,
\begin{equation}\label{eq:gamma-coefficient}
 [z^n](1-z)^\alpha
 =\frac{\Gamma(n-\alpha)}{\Gamma(-\alpha)\Gamma(n+1)},
\end{equation}
with the usual interpretation by analytic continuation at nonnegative integral $\alpha$. Differentiating with respect to $\alpha$ gives the coefficients of the corresponding logarithmic terms. The ratio $\Gamma(n-\alpha)/\Gamma(n+1)$ is asymptotic to $n^{-\alpha-1}$. In particular, a term $(1-z)^r(\log(1-z))^\ell$ with integral $r\geq0$ and $\ell\geq1$ has coefficients of order
\begin{equation}\label{eq:logpower-order}
 O\bigl(n^{-r-1}(1+\log n)^{\ell-1}\bigr).
\end{equation}
At a point $\xi$ the same calculation introduces the factor $\xi^{-n}$. These are the singularity-analysis steps in the hybrid argument; Lemma~\ref{lem:darboux} controls the portion that has only boundary smoothness.

\subsubsection{Factoring the infinite product}

To apply the hybrid method, use the polylogarithmic expansion in \eqref{eq:F-polylogs}. The higher the index of a polylogarithm, the more boundary derivatives it possesses. We exploit this by keeping a finite initial part and estimating the remaining tail.

For an integer $J\geq1$, put
\[
 P_J(z)=\exp\left(\sum_{j=1}^Jb_j\Li_{2j}(z^j)\right),
 \qquad
 Q_J(z)=\exp\left(\sum_{j>J}b_j\Li_{2j}(z^j)\right).
\]
Then $F=P_JQ_J$. For $0\leq r\leq2J$, termwise differentiation gives the uniform bound
\[
 \left|\frac{\dd^r}{\dd z^r}\Li_{2j}(z^j)\right|
 \leq j^r\zeta(2j-r),\qquad |z|\leq1,
 \quad j>J.
\]
Terms of degree less than $r$ simply vanish when differentiated. The smallest zeta argument in this estimate is $2$. Equation~\eqref{eq:bj-bound} therefore proves uniform convergence of every differentiated tail through order $2J$. Consequently $Q_J$ is $C^{2J}$ on the closed disk.

Each polylogarithm in $P_J$ continues across the unit circle except at the roots of $z^j=1$. Near those finitely many roots it has the standard logarithmic expansion. For an integer $s\geq2$, this follows, for example, from
\begin{equation}\label{eq:polylog-local}
 \Li_s(e^{-\tau})
 =\frac{(-\tau)^{s-1}}{(s-1)!}
       \bigl(H_{s-1}-\log\tau\bigr)
   +\sum_{\substack{r\geq0\\r\ne s-1}}
       \zeta(s-r)\frac{(-\tau)^r}{r!},
 \qquad |\tau|<2\pi,
\end{equation}
where $H_r=\sum_{j=1}^rj^{-1}$ and the branch is fixed in the indented neighborhood. Its only nonanalytic term has positive order $s-1$. In particular, $P_J$ is bounded near every exceptional point and admits the expansions required in Lemma~\ref{lem:hybrid}. Taking $J=2$ proves that the coefficient error can be made $o(n^{-2})$; taking $J=4$ makes it $o(n^{-4})$. Notice that retaining only the first polylogarithm leaves a $C^2$ factor and, under this version of the hybrid estimate, gives only an $o(n^{-1})$ error. That precision would not justify a main term of order $n^{-2}$.

\subsubsection{The singular expansion and its coefficients}

Set
\[
 R(z)=\sum_{j\geq2}b_j\Li_{2j}(z^j),\qquad t=1-z.
\]
Then $R'(1)=h$, and the first two derivatives at $1$ exist. Formula~\eqref{eq:polylog-local}, together with the smoothness of the remaining tail, gives
\begin{align*}
 \Li_2(z)&=\zeta(2)+t\log t-t
              +\frac{t^2}{2}\log t-\frac{t^2}{4}
              +O(t^3\log t),\\
 R(z)&=R(1)-ht+\frac12R''(1)t^2+O(t^3\log t).
\end{align*}
The second formula may first be read radially; the hybrid factorization just established is what permits its use in coefficient extraction. Since $\log F=\Li_2/2+R$, writing $A=-1/2-h$ and exponentiating yields
\begin{equation}\label{eq:F-local}
 \frac{F(z)}{C}
 =1+\frac12t\log t+At
   +\frac18t^2(\log t)^2
   +\left(\frac14+\frac A2\right)t^2\log t
   +Dt^2+O(t^3(\log t)^3),
\end{equation}
where the value of the analytic coefficient $D$ is not needed.

We specify the error argument because the product need not have an indented continuation. Use $P_4Q_4$ and expand $P_4$ to an order greater than $4$ at all the relevant roots of unity. Multiply by the Taylor jets of $Q_4$ as in Lemma~\ref{lem:hybrid}. At $z=1$, the terms not displayed in \eqref{eq:F-local} that can affect an error of order $n^{-4}(\log n)^2$ start with $t^3(\log t)^r$, $r\leq3$. Terms of degree $4$ have logarithmic degree at most $4$ and give a smaller error by \eqref{eq:logpower-order}. At $z=-1$, the first singular contribution is from $\Li_4(z^2)$; it has order $(1+z)^3\log(1+z)$ and hence coefficients $O(n^{-4})$. At the other exceptional roots, the first possible singular powers have degree at least $5$ and are smaller still. The hybrid remainder is $o(n^{-4})$. Thus transferring the displayed nonanalytic terms in \eqref{eq:F-local} incurs an overall error $O(n^{-4}(\log n)^2)$.

For $n>2$, differentiation of \eqref{eq:gamma-coefficient} gives
\begin{align*}
 [z^n](1-z)\log(1-z)
   &=\frac1{n(n-1)},\\
 [z^n](1-z)^2\log(1-z)
   &=-\frac2{n(n-1)(n-2)},\\
 [z^n](1-z)^2(\log(1-z))^2
   &=\frac{4H_{n-3}-6}{n(n-1)(n-2)}.
\end{align*}
The polynomial terms in \eqref{eq:F-local} have zero coefficients once $n>2$. We obtain
\[
 a_n=C\left\{\frac1{2n(n-1)}
  +\frac{\tfrac12H_{n-3}-\tfrac34+h}
         {n(n-1)(n-2)}\right\}
  +O\left(\frac{(\log n)^2}{n^4}\right).
\]
Finally, $H_{n-3}=\log n+\gamma+O(n^{-1})$ and expansion of the denominators gives \eqref{eq:an-two-terms}. This proves Theorem~\ref{thm:enumeration}.

\section{The Shepp--Lloyd model}\label{sec:shepp}
For completeness, we give the probabilistic derivation that underlies the coefficient manipulations. Fix $0<\rho<1$, and let $P_k^{(\rho)}$ be independent Poisson variables of means $\rho^k/k$. The sum of the means is $-\log(1-\rho)<\infty$, so only finitely many $P_k^{(\rho)}$ are nonzero almost surely. For a sequence $\chi=(c_k)_{k\geq1}$ of nonnegative integers satisfying $\sum_{k\geq1}kc_k=d$,
\begin{equation}\label{eq:shepp}
\P_\rho((P_k^{(\rho)})=\chi)
=(1-\rho)\rho^d\prod_{k\ge1}\frac1{k^{c_k}c_k!}.
\end{equation}
Summing over all such sequences shows that $\sum kP_k^{(\rho)}$ is geometric, with mass $(1-\rho)\rho^d$. Conditional on $\sum kP_k^{(\rho)}=d$, these variables have the same distribution as the cycle counts of a permutation chosen uniformly from $S_d$. This is the Shepp--Lloyd model \cite{SheppLloyd}. For every nonnegative function $\Phi$ of the cycle counts, Tonelli's theorem gives
\begin{equation}\label{eq:sheppidentity}
\frac{1}{1-\rho}\E_\rho\Phi((P_k^{(\rho)}))
=\sum_{d\ge0}\rho^d\E_{S_d}\Phi((c_k)).
\end{equation}
For signed or complex $\Phi$, the same identity holds whenever the absolute expectations are summable.

Condition each $P_k^{(\rho)}$ to be even. At $\rho<1$ the simultaneous event has positive probability $(1-\rho)E(\rho)$, and the conditioned coordinates remain independent. Writing each even count as twice a nonnegative integer and putting $r=\rho^2$, the resulting sequence has distribution
\begin{equation}\label{eq:Boltzmann}
\P_r(\omega)=\frac{r^{|\omega|}w(\omega)}{F(r)}.
\end{equation}
Conditioning $|\omega|=n$ gives exactly \eqref{eq:profilelaw}, independently of $r$. Unlike the unrestricted Poisson family, the conditioned family has a finite normalizing constant at $r=1$:
\[
C=F(1)=\prod_{k\ge1}\cosh(1/k)<\infty,
\]
because $0\le\log\cosh x\le x^2/2$ for real $x$. At this boundary the paired counts are independent variables $M_k$ with
\begin{equation}\label{eq:Mlaw}
 \P(M_k=m)=\frac1{\cosh(1/k)\,k^{2m}(2m)!},\qquad m\geq0.
\end{equation}
Thus $2M_k$ is Poisson with mean $1/k$, conditioned to be even. We do not obtain this law by conditioning independent Poisson$(1/k)$ variables on all of them being even: that infinite event has probability zero, since its probability is
$\lim_{b\to\infty}\exp(-\sum_{k\le b}1/k)\prod_{k\le b}\cosh(1/k)=0$.

\section{Extremal values}\label{sec:extrema}
For cycle counts described by $\omega=(m_k)\in\mathcal C$, the factorization formula reads
\begin{equation}\label{eq:Hprofile}
 H(\omega)=\prod_{k\geq1}k^{m_k}(2m_k-1)!!.
\end{equation}
The sequence records only how many cycles of each length occur; it does not specify which cycles are paired together. The double factorial counts the possible pairings. Put $X_n=\invol^*(\sigma_{2n})$, so that $X_n=H((M_{k,n})_{k\geq1})$.

\begin{proposition}\label{prop:extrema}
For $n\ge1$ the maximum of $\invol^*$ on $\mathcal E_{2n}$ is $(2n-1)!!$, attained only at the identity. For $n>2$ its minimum is $n-1$, attained precisely at cycle type $1^2(n-1)^2$. The minima for $n=1,2$ are $1,2$, respectively.
\end{proposition}
\begin{proof}
Fix $\sigma$. Choosing the first factor $\tau_1$ determines the second as $\sigma\tau_1$, so there are at most $(2n-1)!!$ factorizations. The identity attains this bound. If $\sigma(x)=y\ne x$, choose a perfect matching $\tau$ containing the pair $\{x,y\}$. Then $\sigma\tau$ fixes $y$. This matching cannot be the first factor of a fixed-point-free factorization, proving strict inequality for every nonidentity $\sigma$.

For the minimum, set $f_k(r)=k^r(2r-1)!!$. If $k\ge2$ and $r\ge1$, then $f_k(r)\ge kr$, with equality exactly when $r=1$. Indeed, the inequality is immediate for $r=1$, and $k^{r-1}(2r-1)!!>r$ for $r\ge2$. If at least two distinct lengths $k\ge2$ are occupied, their quantities $km_k$ are integers at least two and cannot both equal two. For integers $a,b\ge2$, $ab\ge a+b$, with equality only for $a=b=2$; successive application shows that the product of these quantities is strictly greater than their sum. Consequently, writing $s=\sum_{k\ge2}km_k$, the product over $k\ge2$ is at least $s$ when $s>0$, with equality precisely when a single such length is occupied once.

If $m_1=0$, we obtain $H(\omega)\ge n$. If $m_1=1$ and $n>2$, then $s=n-1>0$ and $H(\omega)\ge n-1$, with equality exactly when $m_{n-1}=1$ and all other coordinates vanish. Finally, if $m_1=r\ge2$, then $(2r-1)!!\ge r+1$. When $s=0$ this already exceeds $n-1$; when $s>0$ we have $s\ge2$ and
\[
H(\omega)\ge(r+1)s>r+s-1=n-1.
\]
The two remaining small cases follow directly from \eqref{eq:formula}.
\end{proof}

\section{Two long cycles and the factorization limit}\label{sec:limits}

The factorization formula is multiplicative across cycle lengths. To find its
limiting law, we therefore need to isolate the cycles responsible for its
large deterministic factor. Almost all the letters belong to two cycles
of the same length. Their common length supplies a factor asymptotic to
$n$, and the counts of the other cycles supply the random limiting factor.

Recall that $\omega=(m_k)\in\mathcal C$ describes $2m_k$ cycles
of length $k$, and that these cycles contain $2|\omega|$ letters.
For $|\omega|=n$, their probability is $w(\omega)/a_n$.
Theorem~\ref{thm:enumeration} gives $a_n\sim C/(2n^2)$.
The independent variables $M_k$ from \eqref{eq:Mlaw} will give the
limiting numbers of cycles of each remaining length, divided by two.

The sequence $\Omega=(M_k)_{k\ge1}$ belongs to $\mathcal C$ with probability one.
Indeed, for $x\ge0$,
\[
 1-\operatorname{sech}x
 =\int_0^x\operatorname{sech}t\tanh t\dd t
 \le\frac{x^2}{2},
\]
so $\sum_k\P(M_k>0)<\infty$. The first Borel--Cantelli lemma applies. The probability of each possible sequence is
\begin{equation}\label{eq:residual-mass}
 \P(\Omega=\omega)=\frac{w(\omega)}{C}.
\end{equation}
One can also verify normalization directly: monotone convergence of the
finite products gives $\sum_{\omega\in\mathcal C}w(\omega)=C$.

We shall use the almost surely finite random variables
\begin{equation}\label{eq:residual-statistics}
 S=\sum_{k\ge1}kM_k,\qquad
 W=H(\Omega)=\prod_{k\ge1}(2M_k-1)!!\,k^{M_k}.
\end{equation}

Let
\[
 L_n=\max\{k:c_k(\sigma_{2n})>0\}
\]
be the length of a longest cycle. The removal rule applies exactly when
$L_n>n/2$. In that case, exactly two cycles have length $L_n$: every
cycle multiplicity is even, and four cycles of length greater than $n/2$
would require more than $2n$ letters. The same observation rules out
another cycle length greater than $n/2$. Thus these two equal-length
cycles together contain more than half of the $2n$ letters. We remove
both cycles of length $L_n$ and define
\[
 \Omega_n=\left(M_{k,n}-\mathbf1_{\{k=L_n\}}\right)_{k\geq1},
\]
where the indicator is $1$ for $k=L_n$ and $0$ otherwise. Thus the
$k$th entry of $\Omega_n$ is half the number of $k$-cycles that remain.
If $L_n\leq n/2$, we remove no cycles and put $\Omega_n=\ast$.
The inequality is strict: cycles of length exactly $n/2$ are not removed.
Here $\ast$ records that every cycle has length at most $n/2$; it is not a sequence
of cycle counts. In particular, it differs from the all-zero sequence,
which occurs when $\sigma_{2n}$ consists of two $n$-cycles and their removal
leaves no cycles. The possible values of $\Omega_n$ therefore belong to
$\mathcal C\cup\{\ast\}$, whereas $\Omega$ always belongs to $\mathcal C$.

For example, suppose $n=10$. If $\sigma_{20}$ has cycle type $7^23^2$,
removing the two $7$-cycles leaves two $3$-cycles, so
$\Omega_{10}=(0,0,1,0,\ldots)$. Cycle type $10^2$ gives the all-zero
sequence after removal. In contrast, cycle type $5^4$ gives
$\Omega_{10}=\ast$, since no cycle length exceeds $5$.

Write $\mathcal L(Z)$ for the distribution of a random variable $Z$.
If two probability distributions $\mu$ and $\nu$ assign probabilities to
the same countable set of possible values, their total variation distance is
\[
 d_{\TV}(\mu,\nu)=\frac12\sum_x|\mu\{x\}-\nu\{x\}|,
\]
where the sum is over all those values. For $\Omega_n$ and $\Omega$, we
sum over $\mathcal C\cup\{\ast\}$ and take $\P(\Omega=\ast)=0$. Thus
the comparison includes the probability that no cycle length exceeds $n/2$.

\begin{theorem}\label{thm:residual}
The distribution of $\Omega_n$, defined above by removing both cycles of length $L_n$ when $L_n>n/2$, converges in total variation to the distribution of $\Omega$.
Consequently,
\[
 \P(L_n>n/2)\longrightarrow1,
 \qquad n-L_n\ \xrightarrow{\mathrm d}\ S,
 \qquad \frac{L_n}{n}\ \xrightarrow{\P}\ 1.
\]
\end{theorem}

\begin{proof}
Fix $\omega=(m_k)\in\mathcal C$ and put $s=|\omega|$. If $n>2s$,
the event $\Omega_n=\omega$ says exactly that $\sigma_{2n}$ has $2m_k$
cycles of length $k$ for each $k$, together with two additional
$(n-s)$-cycles. The length $n-s$ exceeds both $n/2$ and every $k$
for which $m_k>0$. Formula~\eqref{eq:profilelaw} therefore
gives the exact identity
\begin{equation}\label{eq:residual-exact}
 \P(\Omega_n=\omega)
 =\frac{w(\omega)}{2(n-s)^2a_n}.
\end{equation}
By Theorem~\ref{thm:enumeration}, this tends to $w(\omega)/C$.

For completeness, the passage from pointwise convergence to total variation
requires no uniform estimate in $s$. Given $\varepsilon>0$, choose a finite
subset $\mathcal B$ of $\mathcal C$ such that
\[
 \sum_{\omega\in\mathcal B}\frac{w(\omega)}{C}>1-\varepsilon.
\]
Pointwise convergence on $\mathcal B$ implies
$\P(\Omega_n\in\mathcal B)>1-2\varepsilon$ for all sufficiently large $n$.
The sum of the absolute differences of the masses on $\mathcal B$ tends
to zero, and their sum on the complement, including $\ast$, is at most
$3\varepsilon$ for all sufficiently large $n$. Letting $\varepsilon$ tend
to zero proves total variation convergence. In particular,
$\P(\Omega_n=\ast)\to0$.

When $L_n>n/2$, the cycles left after removal contain $2(n-L_n)$
letters, so $n-L_n=|\Omega_n|$. The preceding convergence gives
$n-L_n\xrightarrow{\mathrm d}S$. Since $S$ is finite with probability
one, $n-L_n$ is bounded in probability, and hence $L_n/n\to1$ in probability.
\end{proof}

\subsection{The limiting factorization count}

\begin{theorem}\label{thm:factorization-limit}
With $W$ as in \eqref{eq:residual-statistics},
\begin{equation}\label{eq:factorization-limit}
 \frac{X_n}{n}\xrightarrow{\mathrm d}W.
\end{equation}
For each positive integer $j$, the probability that $W=j$ is
\begin{equation}\label{eq:W-pmf}
 \P(W=j)=\frac1C
 \sum_{\substack{\omega\in\mathcal C\\H(\omega)=j}}w(\omega).
\end{equation}
The sum in \eqref{eq:W-pmf} is finite. In particular,
\begin{equation}\label{eq:limit-atoms}
 \P(W=1)=\frac{3}{2C}.
\end{equation}
\end{theorem}
\begin{proof}
When $L_n>n/2$, exactly two cycles have length $L_n$, and every
cycle left after their removal has smaller length. The factorization formula gives
\begin{equation}\label{eq:statistics-residual}
 X_n=(n-|\Omega_n|)H(\Omega_n).
\end{equation}
Choose a finite subset $\mathcal B$ of $\mathcal C$ for which
$\P(\Omega\in\mathcal B)>1-\varepsilon$. On this set,
\[
 (1-|\omega|/n)H(\omega)\longrightarrow H(\omega)
\]
uniformly. The total variation convergence of $\Omega_n$ in
Theorem~\ref{thm:residual}, followed by bounded continuous test functions,
proves \eqref{eq:factorization-limit}. The complement of $\mathcal B$,
including $\ast$, contributes at most a constant times
$\varepsilon+o(1)$, and $\varepsilon$ is arbitrary.

Summing \eqref{eq:residual-mass} over $H(\omega)=j$ proves
\eqref{eq:W-pmf}. Every occupied index satisfies $k\le j$, while
$(2m_k-1)!!\le j$ bounds its multiplicity. There are consequently only
finitely many contributing sequences. Finally, $H(\omega)=1$ holds
only for the all-zero sequence and the sequence $m_1=1$, $m_k=0$ for $k\geq2$. Their weights are
$1$ and $1/2$, which proves \eqref{eq:limit-atoms}.
\end{proof}

\subsection{The rate of convergence of the remaining cycle counts}

The exact probability in \eqref{eq:residual-exact} also determines
the rate of convergence. There are two sources of error. First, some
sequences describe too many letters to remain after two cycles of length
greater than $n/2$ are removed. Second, the probabilities of the sequences
that can remain differ from their limiting probabilities. This second
error is larger. The difference changes sign when the number of remaining
letters is of order $\log n$, producing a second logarithm in the
total variation estimate.

\begin{theorem}\label{thm:residual-rate}
For the variables $\Omega_n$ and $\Omega$ defined above,
\[
 d_{\TV}\bigl(\mathcal L(\Omega_n),\mathcal L(\Omega)\bigr)
 =\frac{\log n-\log\log n+\log2-\tfrac12+o(1)}{n}.
\]
In particular, the leading equivalent is $(\log n)/n$.
\end{theorem}
\begin{proof}
Put $q_s=\P(S=s)=a_s/C$, including $q_0=1/C$, and define
\[
 \alpha_n=\frac{C}{2n^2a_n},\qquad
 \theta_n=n(1-\sqrt{\alpha_n}),\qquad
 \Lambda_n=\log n+\kappa.
\]
The enumeration theorem gives
\begin{equation}\label{eq:rate-normalizer}
 \alpha_n=1-\frac{\Lambda_n}{n}
             +O\left(\frac{\log^2n}{n^2}\right),\qquad
 \theta_n=\frac{\Lambda_n}{2}
             +O\left(\frac{\log^2n}{n}\right).
\end{equation}
In particular, $0<\theta_n<n/2$ for all sufficiently large $n$.

For $\omega\in\mathcal C$ with $s=|\omega|<n/2$, the probability
that these counts remain after removal, divided by its limiting value, is
\[
 \frac{\P(\Omega_n=\omega)}{\P(\Omega=\omega)}
 =\frac{\alpha_n}{(1-s/n)^2}.
\]
This ratio is increasing in $s$ and is less than $1$ precisely when
$s<\theta_n$. A sequence with $|\omega|\geq n/2$ cannot be left
after removal of the two long cycles, so $\P(\Omega_n=\omega)=0$.
Also $\P(\Omega=\ast)=0$. Total variation equals the sum of the positive
differences $\P(\Omega=\omega)-\P(\Omega_n=\omega)$; the symbol $\ast$
contributes no positive difference. We therefore obtain the exact identity
\begin{equation}\label{eq:tv-exact-threshold}
 d_{\TV}\bigl(\mathcal L(\Omega_n),\mathcal L(\Omega)\bigr)
 =\P(S\geq n/2)
  +\sum_{0\leq s<\theta_n}q_s
       \left(1-\frac{\alpha_n}{(1-s/n)^2}\right).
\end{equation}
This identity does not require an estimate of $\P(L_n\leq n/2)$.

We record the two truncated sums needed to evaluate it. The coefficient
expansion implies
\[
 q_s=\frac1{2s^2}+O\left(\frac{\log(s+1)}{s^3}\right)
 \qquad(s\geq1).
\]
Consequently, as a real parameter $t$ tends to infinity,
\begin{align}
 \P(S\geq t)&=\frac1{2t}
              +O\left(\frac{\log t}{t^2}\right),
                 \label{eq:residual-tail-rate}\\
 \sum_{0\leq s<t}s q_s
   &=\frac12\log t+\frac\gamma2+h
              +O\left(\frac{\log t}{t}\right).
                 \label{eq:residual-truncated-mean}
\end{align}
To identify the constant in the second line, note that
$\sum_{s\geq1}|s q_s-1/(2s)|<\infty$. Also
\[
 \frac{F'(r)}C
 =\frac{F(r)}C\left(-\frac{\log(1-r)}{2r}+R'(r)\right),
 \qquad 0<r<1.
\]
The local expansion of $F$ gives
$(F(r)/C-1)\log(1-r)\to0$, while $R'(r)\to h$. Taking the
limit in the absolutely convergent difference of series yields
\[
 \sum_{s\geq1}\left(sq_s-\frac1{2s}\right)=h.
\]
The harmonic-number expansion and the bound on the tail of this
absolutely convergent series prove
\eqref{eq:residual-truncated-mean}. Both estimates are valid with the
strict cutoff $s<t$; integer rounding contributes only their stated
error terms.

Uniformly for $0\leq s<\theta_n=O(\log n)$,
\[
 1-\frac{\alpha_n}{(1-s/n)^2}
 =\frac{\Lambda_n-2s}{n}
   +O\left(\frac{\log^2n}{n^2}\right).
\]
Thus \eqref{eq:tv-exact-threshold} and the two truncated estimates give
\begin{align*}
 n\,d_{\TV}\bigl(\mathcal L(\Omega_n),\mathcal L(\Omega)\bigr)
 &=1+\Lambda_n\left(1-\frac1{2\theta_n}
          +O\left(\frac{\log\theta_n}{\theta_n^2}\right)\right)\\
 &\qquad-\log\theta_n-\gamma-2h
       +O\left(\frac{\log\theta_n}{\theta_n}
                +\frac{\log^2n}{n}\right)\\
 &=\Lambda_n-\log\theta_n-\gamma-2h+o(1).
\end{align*}
Here the initial $1$ cancels $\Lambda_n/(2\theta_n)=1+o(1)$.
Finally,
$\log\theta_n=\log\log n-\log2+o(1)$ and
$\kappa=\gamma-\tfrac12+2h$. Substitution proves the theorem.
In particular, the constant does not require a new infinite product
or series: the constants already present in the two-term enumeration
cancel.
\end{proof}

\begin{proposition}\label{prop:no-giant-and-count}
The probability that no cycle length exceeds $n/2$ satisfies
\begin{equation}\label{eq:no-giant-rate}
 \P(L_n\leq n/2)=O\left(\frac{\log^2n}{n^2}\right).
\end{equation}
Using the variables $S$ and $W$ from \eqref{eq:residual-statistics}, define
\[
 \widehat X_n=
 \begin{cases}
 (n-S)W,&S<n/2,\\
 0,&S\geq n/2.
 \end{cases}
\]
Then
\begin{equation}\label{eq:count-TV-approximation}
 d_{\TV}\bigl(\mathcal L(X_n),\mathcal L(\widehat X_n)\bigr)
 \leq\frac{\log n-\log\log n+\log2-\tfrac12+o(1)}n.
\end{equation}
\end{proposition}
\begin{proof}
Retain $q_s=a_s/C$ and $\alpha_n=C/(2n^2a_n)$ from the preceding
proof. Summing the exact probabilities in \eqref{eq:residual-exact} gives
\[
 \P(L_n\leq n/2)
 =1-\alpha_n\sum_{0\leq s<n/2}\frac{q_s}{(1-s/n)^2}.
\]
To estimate this sum, put $J(x)=(1-x)^{-2}-1-2x$. The quotient
$J(x)/x^2$ extends smoothly to $[0,1/2]$ with value $3$ at zero.
Since $q_s=1/(2s^2)+O(\log(s+1)/s^3)$, a Riemann sum gives
\begin{align*}
 \sum_{1\leq s<n/2}q_sJ(s/n)
 &=\frac1{2n}\int_0^{1/2}\frac{J(x)}{x^2}\,\dd x
       +O\left(\frac{\log^2n}{n^2}\right)\\
 &=\frac{\tfrac12+\log2}{n}
       +O\left(\frac{\log^2n}{n^2}\right).
\end{align*}
Indeed, replacing $q_s$ by $1/(2s^2)$ incurs an error bounded by
$(K/n^2)\sum_{s<n/2}\log(s+1)/s=O(n^{-2}\log^2n)$ for a constant
$K$. The smooth Riemann sum has error $O(n^{-2})$. The integral
follows from $J(x)/x^2=2/(1-x)+(1-x)^{-2}$.
Together with \eqref{eq:residual-tail-rate} and
\eqref{eq:residual-truncated-mean}, this yields
\begin{align*}
 \sum_{0\leq s<n/2}\frac{q_s}{(1-s/n)^2}
 &=\P(S<n/2)+\frac2n\E(S;S<n/2)
       +\sum_{1\leq s<n/2}q_sJ(s/n)\\
 &=1+\frac{\log n+\gamma+2h-\tfrac12}{n}
       +O\left(\frac{\log^2n}{n^2}\right)\\
 &=1+\frac{\log n+\kappa}{n}
       +O\left(\frac{\log^2n}{n^2}\right).
\end{align*}
Its first correction cancels that of $\alpha_n$ in
\eqref{eq:rate-normalizer}, proving \eqref{eq:no-giant-rate}.

Finally, assign to each $\omega\in\mathcal C$ with $s=|\omega|<n/2$
the number $(n-s)H(\omega)$. Assign zero to $\ast$ and to every sequence
with $|\omega|\geq n/2$. This rule gives $\widehat X_n$ when applied
to $\Omega$, and agrees with $X_n$ when applied to $\Omega_n$ on the
event $L_n>n/2$. Total
variation contracts under a map, while altering a random variable
on an event of probability $\varepsilon$ changes its law by at
most $\varepsilon$ in total variation. Thus the left side of
\eqref{eq:count-TV-approximation} is at most the total
variation distance between $\Omega_n$ and $\Omega$ plus $\P(L_n\leq n/2)$. The latter is $o(n^{-1})$,
so Theorem~\ref{thm:residual-rate} proves the assertion.
\end{proof}

\begin{remark}
The convergence in \eqref{eq:factorization-limit} does not imply uniform
convergence of the distribution functions, or convergence in total
variation. For $n>2$, the cycle type $1^2(n-1)^2$ has $X_n=n-1$ and
probability
\[
 \frac{1}{4(n-1)^2a_n}\longrightarrow\frac1{2C}.
\]
Since $W\geq1$, we have
\[
 \P\left(W\leq1-\frac1{2n}\right)=0,
 \qquad
 \P\left(\frac{X_n}{n}\leq1-\frac1{2n}\right)
 \geq\frac{1}{4(n-1)^2a_n}.
\]
Thus the maximum difference between the two distribution functions
stays bounded away from zero. Total variation is at least this
difference, so it cannot tend to zero either.
\end{remark}

\subsection{The median number of factorizations}

The median can be determined at every degree. There is one exceptional
half-degree, namely $n=3$, at which the median is $2$. To establish that
there are no later exceptions, we need a uniform upper bound for the
normalizing coefficient $a_n$.

\begin{lemma}\label{lem:global-density}
For every integer $n\ge1$,
\[
 a_n<\frac{3}{2n^2}.
\]
\end{lemma}
\begin{proof}
Write $f\preceq g$ when every coefficient of $g-f$ is nonnegative. Since
$(2m)!\ge2^m m!$ for every $m\ge0$, the product over cycle lengths gives
\[
 F(z)\preceq B(z):=\exp\left(\frac12\Li_2(z)\right)
       =\sum_{n\ge0}\beta_nz^n.
\]
We shall prove $n^2\beta_n<3/2$. Logarithmic differentiation gives the
positive recurrence
\begin{equation}\label{eq:majorant-recurrence}
 \beta_0=1,\qquad
 n\beta_n=\frac12\sum_{j=1}^n\frac{\beta_{n-j}}j.
\end{equation}
Also,
\begin{equation}\label{eq:majorant-total}
 \sum_{n\ge0}\beta_n=B(1)=\exp(\pi^2/12)<\frac73.
\end{equation}
For an elementary numerical bound in the last step, use $\pi^2<10$ and
bound the tail of the exponential series after its cubic term by a
geometric series:
\[
 \ee^{5/6}
 <1+\frac56+\frac{25}{72}+\frac{125}{1296}
       +\frac{625}{31104}\frac1{1-1/6}
 =\frac{11929}{5184}<\frac73.
\]

The first fourteen cases of the induction admit a short integer
certificate. Define $u_0=10000$ and, for $1\le n\le14$, put
\[
 u_n=\left\lceil\frac1{2n}\sum_{j=1}^n\frac{u_{n-j}}j\right\rceil.
\]
The positivity of \eqref{eq:majorant-recurrence} implies
$\beta_n\le u_n/10000$. The exact integer values are
\[
\begin{aligned}
 (u_1,\ldots,u_7)&=(5000,2500,1389,851,563,396,291),\\
 (u_8,\ldots,u_{14})&=(222,175,141,116,96,82,70).
\end{aligned}
\]
For each of these indices, $n^2u_n\le14259<15000$, proving the required
base cases without a floating-point comparison.

Suppose now that $n\ge15$ and $k^2\beta_k<3/2$ for $1\le k<n$.
Reindexing \eqref{eq:majorant-recurrence} and writing
$n/(n-k)=1+k/(n-k)$ gives
\begin{align*}
 n^2\beta_n
 &=\frac12\left(1+\sum_{k=1}^{n-1}\beta_k
               +\sum_{k=1}^{n-1}\frac{k\beta_k}{n-k}\right)\\
 &<\frac76+\frac34\sum_{k=1}^{n-1}\frac1{k(n-k)}
 =\frac76+\frac{3H_{n-1}}{2n}.
\end{align*}
The sequence $H_{n-1}/n$ decreases for $n\ge3$, and
\[
 \frac{H_{14}}{15}=\frac{1171733}{5405400}<\frac29.
\]
Consequently $n^2\beta_n<7/6+1/3=3/2$. Induction and
$a_n\le\beta_n$ complete the proof.
\end{proof}

\begin{proposition}\label{prop:low-factorization-counts}
For every $n\ge3$,
\begin{align}
 \P(X_n<n)&=\P(X_n=n-1)
       =\frac1{4(n-1)^2a_n},\label{eq:exact-lower-median-tail}\\
 \P(X_n=n)&=\frac1{2n^2a_n}
       +\frac{\mathbf1_{\{n=6\}}}{288a_n}.
       \label{eq:exact-median-atom}
\end{align}
The first event consists of cycle type $1^2(n-1)^2$. The second consists
of cycle type $n^2$, with the additional type $1^22^23^2$ when $n=6$.
\end{proposition}
\begin{proof}
The first assertion, including the uniqueness of the minimizing type,
follows from Proposition~\ref{prop:extrema}. The weight of these cycle counts is
$1/[4(n-1)^2]$.

For the second assertion, write $r=m_1$ and
$s=\sum_{k\ge2}km_k$. If $r=0$, the proof of
Proposition~\ref{prop:extrema} shows that $H(\omega)\ge n$, with
equality precisely when $m_n=1$. If $r\ge2$, then
$(2r-1)!!\ge r+1$. When $s=0$ this exceeds $n=r$; when $s>0$, we have
$s\ge2$ and
\[
 H(\omega)\ge(r+1)s>r+s=n.
\]
It remains to consider $r=1$, so that $s=n-1$.

List the occupied lengths at least two as $k_1,\ldots,k_t$, and put
$v_i=k_im_{k_i}$. Then $v_i\ge2$ and
$H(\omega)\ge\prod_i v_i$. If $t\ge3$, the product of these
integers exceeds their sum by at least two: for three factors this
follows from $abc-a-b-c\ge2$ for $a,b,c\ge2$, and adjoining further
factors preserves the inequality. Thus $H(\omega)>s+1=n$.

If $t=1$ and $m_{k_1}=1$, the count is $s=n-1$. If $t=1$ and
$m_{k_1}\ge2$, the elementary inequality
\[
 k^m(2m-1)!!\ge3km\qquad(k\ge2,\ m\ge2)
\]
shows that the count exceeds $s+1$. The inequality holds at $m=2$ and
then follows by multiplication by $k(2m+1)$.

Finally suppose $t=2$. If a multiplicity is at least two, the
corresponding $v_i$ is at least four, and
$v_1v_2-v_1-v_2\ge2$. Otherwise the two multiplicities are both one,
and the equality $H(\omega)=s+1$ is
\[
 k_1k_2=k_1+k_2+1,
 \qquad\text{or equivalently}\qquad
 (k_1-1)(k_2-1)=2.
\]
Since the lengths are distinct and at least two, their set must be
$\{2,3\}$. Together with $m_1=1$, this gives $n=6$. Its weight is
$(2!\,2^2 2!\,3^2 2!)^{-1}=1/288$, whereas the weight of $m_n=1$
is $1/(2n^2)$. Division by $a_n$ proves
\eqref{eq:exact-median-atom}.
\end{proof}

\begin{theorem}\label{thm:sharp-median}\label{thm:median}
The unique median of $X_n$ is $n$ for every $n\ge4$. The unique medians
at $n=1,2,3$ are $1,2,2$, respectively, so the threshold $4$ is sharp.
Moreover,
\begin{equation}\label{eq:median-cdf}
 \P(X_n<n)\longrightarrow\frac1{2C},\qquad
 \P(X_n\le n)\longrightarrow\frac3{2C},\qquad
 \P(X_n=n)\longrightarrow\frac1C.
\end{equation}
\end{theorem}
\begin{proof}
For $n\ge4$, the two distinct sequences $m_n=1$ and
$m_1=m_{n-1}=1$ give
\[
 a_n\ge\frac1{2n^2}+\frac1{4(n-1)^2}
       >\frac1{2(n-1)^2},
\]
because $n^2<2(n-1)^2$. Equation~\eqref{eq:exact-lower-median-tail}
therefore gives $\P(X_n<n)<1/2$. On the other hand,
Lemma~\ref{lem:global-density} gives
\[
 \frac1{2n^2}+\frac1{4(n-1)^2}
 >\frac3{4n^2}>\frac{a_n}{2}.
\]
Both sequences give factorization counts at most $n$, so
$\P(X_n\le n)>1/2$. These strict inequalities make $n$ the unique
median. The case $n=1$ is immediate; the cases $n=2,3$ follow from the exact distributions in Examples~\ref{ex:median-four-letters} and~\ref{ex:median-six-letters} below. Finally, the asserted limits follow directly from
\eqref{eq:exact-lower-median-tail}--\eqref{eq:exact-median-atom} and
$a_n\sim C/(2n^2)$.
\end{proof}

The small cases illustrate the exceptional half-degree and check the weights of the cycle counts directly.

\begin{example}\label{ex:median-four-letters}
When $n=2$, the possible cycle types are $2^2$ and $1^4$.
Their weights are $1/8$ and $1/24$, and their factorization counts are
$2$ and $3$, respectively. Since $a_2=1/6$,
\[
 \P(X_2=2)=\frac34,\qquad \P(X_2=3)=\frac14.
\]
The median is therefore $2$.
\end{example}

\begin{example}\label{ex:median-six-letters}
When $n=3$, the possible cycle types $3^2$, $1^22^2$, and $1^6$ have weights
$1/18$, $1/16$, and $1/720$, respectively. Their factorization counts are
$3$, $2$, and $15$. Thus
\[
 a_3=\frac1{18}+\frac1{16}+\frac1{720}=\frac{43}{360},
\]
and
\[
 \P(X_3=2)=\frac{45}{86},\qquad
 \P(X_3=3)=\frac{20}{43},\qquad
 \P(X_3=15)=\frac1{86}.
\]
Here the median is $2$, rather than $3$.
\end{example}

Numerically, the limiting probabilities in \eqref{eq:median-cdf} are
$0.236243$, $0.708729$, and $0.472486$, respectively. The large atom at
$n$ has a simple source: it is the cycle type with just two $n$-cycles.
The smaller atom at $n-1$ comes from adjoining two fixed points to two
$(n-1)$-cycles.

All these statements concern the uniform distribution on
$\mathcal E_{2n}$. For a uniform permutation of the whole symmetric group
$S_{2n}$, the factorization count vanishes with probability $1-a_n\to1$,
so its unique median is instead zero for all sufficiently large $n$.

\section{The first two moments of the factorization count}\label{sec:moments}

The limiting law describes the factorization count on the scale on which
most permutations occur. Its moments tell a different story. In this
section we obtain exact generating functions for the first two moments
and evaluate their coefficients. The second moment leads to a fractional
power of the partition generating function. Thus its asymptotic analysis
requires an essential-singularity calculation, rather than the
logarithmic singularities used to enumerate acme permutations.

Throughout this section, $X_n=\invol^*(\sigma_{2n})$, and expectations refer
to the uniform distribution on acme permutations of $[2n]$.

\begin{proposition}\label{prop:momentproducts}
For every real $q$, the following identity holds as an identity of formal
power series:
\begin{equation}\label{eq:general-moment-product}
 \sum_{n\geq0}a_n\E X_n^q\,z^n
 =\prod_{k\geq1}\left(
   \sum_{m\geq0}
   \frac{((2m-1)!!)^q}{(2m)!}
   k^{(q-2)m}z^{km}\right),
\end{equation}
where the term for $n=0$ is interpreted as $1$. In particular, put
\begin{equation}\label{eq:moment-numerators}
 d_n=4^{-n}\binom{2n}{n},\qquad
 \mathcal H(z)=\prod_{k\geq1}(1-z^k)^{-1/2},\qquad h_n=[z^n]\mathcal H(z).
\end{equation}
Then
\begin{equation}\label{eq:exact-two-moments}
 \E X_n=\frac{d_n}{a_n},\qquad
 \E X_n^2=\frac{h_n}{a_n}.
\end{equation}
\end{proposition}

\begin{proof}
For $\omega=(m_k)\in\mathcal C$ with $|\omega|=n$, the probability
of these cycle counts is
\[
 \frac{1}{a_n}\prod_{k\geq1}\frac{1}{k^{2m_k}(2m_k)!},
\]
whereas the factorization count is
$\prod_{k\geq1}(2m_k-1)!!\,k^{m_k}$.
Multiplying these expressions and summing over $\omega\in\mathcal C$ with $|\omega|=n$ proves
\eqref{eq:general-moment-product}. Each coefficient involves only
finitely many sequences, so this step requires no assertion about the
analytic convergence of the product for general $q$.

When $q=1$, the identity
$(2m-1)!!/(2m)!=1/(2^m m!)$ gives
\[
 \sum_{n\geq0}a_n\E X_n\,z^n
 =\prod_{k\geq1}\exp\left(\frac{z^k}{2k}\right)
 =(1-z)^{-1/2}.
\]
When $q=2$, we instead have
\[
 \frac{((2m-1)!!)^2}{(2m)!}
 =4^{-m}\binom{2m}{m},
\]
and the binomial series gives the product $\mathcal H(z)$ in
\eqref{eq:moment-numerators}. This proves
\eqref{eq:exact-two-moments}.
\end{proof}

There is also a direct interpretation of the first identity in
\eqref{eq:exact-two-moments}. Summing $\invol^*(\sigma)$ over all
permutations counts every ordered pair of fixed-point-free involutions
exactly once. There are $((2n-1)!!)^2$ such pairs, and
\[
 \frac{((2n-1)!!)^2}{(2n)!}=4^{-n}\binom{2n}{n}=d_n.
\]

We next give the coefficient estimate needed for the second moment. The
local expansion is of the type studied by Wright~\cite{Wright}; see also
the treatment of partition products in~\cite[Chapter~VIII.6]{FS}.
Because $\mathcal H$ has singularities at every root of unity, a local expansion
at $1$ alone does not justify the conclusion. The proof includes a bound
for the remainder of the coefficient contour.

\begin{lemma}\label{lem:partitionhalf}
As $n\to\infty$,
\begin{equation}\label{eq:half-partition-asymptotic}
 h_n\sim\frac{1}{4\,3^{3/8}}\,
 n^{-7/8}\exp\left(\pi\sqrt{\frac n3}\right).
\end{equation}
\end{lemma}

\begin{proof}
Set
\[
 A=\frac{\pi^2}{12},\qquad \lambda=(2\pi)^{-1/4}.
\]
The partition product
$\mathcal P(z)=\prod_{k\geq1}(1-z^k)^{-1}$ satisfies the classical
Dedekind transformation
\begin{equation}\label{eq:partition-transform}
 \mathcal P(\ee^{-t})
 =\left(\frac{t}{2\pi}\right)^{1/2}
   \exp\left(\frac{\pi^2}{6t}-\frac{t}{24}\right)
   \mathcal P(\ee^{-4\pi^2/t}),\qquad \Re t>0;
\end{equation}
see~\cite[Chapter~VIII.6, Eq.~(70)]{FS} or
\cite[\S27.14(iv)]{DLMF}.
Define $\mathcal H$ in the unit disc by exponentiating
$-\frac12\sum_{k\geq1}\log(1-z^k)$, with the logarithms normalized at
$z=0$. It is analytic and nonzero there. Taking the square root of
\eqref{eq:partition-transform}, initially for positive real $t$ and
then by analytic continuation throughout the right half-plane, gives
\begin{equation}\label{eq:half-partition-transform}
 \mathcal H(\ee^{-t})
 =\lambda t^{1/4}\exp\left(\frac At-\frac{t}{48}\right)
   \mathcal H(\ee^{-4\pi^2/t}).
\end{equation}
Here $t^{1/4}$ denotes the principal branch. In particular, uniformly
for $t\to0$ with $|\arg t|\leq\pi/4$,
\begin{equation}\label{eq:half-partition-local}
 \mathcal H(\ee^{-t})
 =\lambda t^{1/4}\exp\left(\frac At-\frac{t}{48}\right)
   \left(1+O(\ee^{-c/|t|})\right)
\end{equation}
for some $c>0$. Indeed,
$\Re(1/t)\geq1/(\sqrt2|t|)$ in this sector, and $\mathcal H(u)=1+O(u)$ at
$u=0$.

Choose
\[
 \alpha=\sqrt{\frac An},\qquad r=\ee^{-\alpha},\qquad
 \theta_0=\alpha^{7/5}.
\]
Cauchy's formula takes the form
\begin{equation}\label{eq:half-partition-cauchy}
 h_n=\frac{1}{2\pi}\int_{-\pi}^{\pi}
 \mathcal H(\ee^{-\alpha+i\theta})
 \exp\bigl(n(\alpha-i\theta)\bigr)\,d\theta.
\end{equation}
First suppose $|\theta|\leq\theta_0$ and put $t=\alpha-i\theta$.
The choice $n=A/\alpha^2$ cancels the linear term in the exponent:
\begin{equation}\label{eq:saddle-expansion}
 \frac At+nt
 =\frac{2A}{\alpha}-\frac{A\theta^2}{\alpha^3}
   +O\left(\frac{|\theta|^3}{\alpha^4}\right).
\end{equation}
The error in~\eqref{eq:saddle-expansion} is $O(\alpha^{1/5})$
uniformly on this arc. Moreover,
\[
 t^{1/4}=\alpha^{1/4}\bigl(1+O(\alpha^{2/5})\bigr),\qquad
 \ee^{-t/48}=1+O(\alpha).
\]
Substituting these estimates into
\eqref{eq:half-partition-local} shows that the contribution of the
central arc to~\eqref{eq:half-partition-cauchy} is
\[
 \frac{\lambda\alpha^{1/4}}{2\pi}\ee^{2A/\alpha}
 \int_{-\theta_0}^{\theta_0}
       \exp\left(-\frac{A\theta^2}{\alpha^3}\right)\,d\theta
 \,(1+o(1)).
\]
Since $\theta_0/\alpha^{3/2}=\alpha^{-1/10}\to\infty$, the
Gaussian integral is asymptotic to
$\sqrt{\pi/A}\,\alpha^{3/2}$. Thus the central contribution is
\begin{equation}\label{eq:half-partition-central}
 \frac{\lambda}{2\sqrt\pi}A^{3/8}n^{-7/8}
 \ee^{2\sqrt{An}}(1+o(1)).
\end{equation}

It remains to bound the rest of the contour. Write
$\sigma_{-1}(j)=\sum_{d\mid j}d^{-1}$. Absolute convergence inside
the unit disc gives
\[
 \log \mathcal H(z)=\frac12\sum_{j\geq1}\sigma_{-1}(j)z^j.
\]
As $\sigma_{-1}(j)\geq1$, we have for $-\pi\leq\theta\leq\pi$
\begin{align}
 \log \mathcal H(r)-\log|\mathcal H(r\ee^{i\theta})|
 &\geq\frac12\sum_{j\geq1}r^j(1-\cos(j\theta))\notag\\
 &=\frac{r(1+r)(1-\cos\theta)}
 {2(1-r)\bigl((1-r)^2+2r(1-\cos\theta)\bigr)}.
 \label{eq:half-partition-modulus}
\end{align}
The last equality follows by summing two geometric series. Since
$1-r\asymp\alpha$ and
$2\theta^2/\pi^2\leq1-\cos\theta\leq\theta^2/2$ on this interval,
there is a constant $c_1>0$ such that
\[
 \log \mathcal H(r)-\log|\mathcal H(r\ee^{i\theta})|
 \geq c_1\frac{\theta^2}{\alpha(\alpha^2+\theta^2)}.
\]
The expression on the right increases with $|\theta|$. Consequently,
for $\theta_0\leq|\theta|\leq\pi$ and all sufficiently small
$\alpha$,
\[
 |\mathcal H(r\ee^{i\theta})|
 \leq \mathcal H(r)\exp(-c_2\alpha^{-1/5})
\]
with $c_2>0$. The absolute value of the contribution of these arcs to
\eqref{eq:half-partition-cauchy} is therefore at most
\[
 \ee^{n\alpha}\mathcal H(r)\exp(-c_2\alpha^{-1/5}).
\]
By~\eqref{eq:half-partition-local}, this is negligible relative to
\eqref{eq:half-partition-central}: the central contribution has the
same exponential factor $\ee^{n\alpha}\mathcal H(r)$, multiplied by a factor
of order $\alpha^{3/2}$. This estimate treats all noncentral arcs
simultaneously, including those containing other roots of unity.

Finally,
\[
 \frac{\lambda}{2\sqrt\pi}A^{3/8}
 =\frac{(2\pi)^{-1/4}}{2\sqrt\pi}
   \left(\frac{\pi^2}{12}\right)^{3/8}
 =\frac{1}{4\,3^{3/8}},
 \qquad 2\sqrt{An}=\pi\sqrt{\frac n3}.
\]
This proves~\eqref{eq:half-partition-asymptotic}.
\end{proof}

\begin{theorem}\label{thm:rawmoments}
For a uniformly chosen acme permutation of $[2n]$,
\begin{align}
 \E X_n
 &\sim\frac{2}{C\sqrt\pi}\,n^{3/2},
 \label{eq:first-moment-asymptotic}\\
 \E X_n^2\sim\Var(X_n)
 &\sim\frac{1}{2C\,3^{3/8}}\,
 n^{9/8}\exp\left(\pi\sqrt{\frac n3}\right).
 \label{eq:second-moment-asymptotic}
\end{align}
\end{theorem}

\begin{proof}
Stirling's formula gives $d_n\sim(\pi n)^{-1/2}$. Combine this with
$a_n\sim C/(2n^2)$ in~\eqref{eq:exact-two-moments} to obtain
\eqref{eq:first-moment-asymptotic}. The second identity in
\eqref{eq:exact-two-moments} and Lemma~\ref{lem:partitionhalf} give
the equivalent for $\E X_n^2$. Finally, $(\E X_n)^2$ has order
$n^3$, which is negligible relative to the right-hand side
of~\eqref{eq:second-moment-asymptotic}. Subtracting it proves the
variance equivalent.
\end{proof}

The scales in Theorem~\ref{thm:rawmoments} explain why moment
calculations alone give a misleading picture of a typical acme
permutation. The factorization count has a discrete limit after division
by $n$, but its expectation is of order $n^{3/2}$, and its second moment
is larger still. Rare sequences of cycle counts carry enough factorizations to dominate
these averages. In particular, the family $(X_n/n)$ is not uniformly
integrable.

\section{Fractional and higher moments}\label{sec:higher}

The first two moments already show that the distributional limit does not
control the mean. We now make this distinction more precise. Moments of order
less than one do pass to the limit. Between orders one and two, an exponential
singularity leads to a saddle-point estimate. Above order two, the identity
permutation alone supplies the leading term. These regimes use different
features of the exact factorization formula.

Write
\[
 r_m=(2m-1)!!,\qquad
 v_{q,m}=\frac{r_m^q}{(2m)!},\qquad r_0=1.
\]
The local series and the moment numerator are
\begin{equation}\label{eq:local-moment-series}
 g_q(u)=\sum_{m\geq0}v_{q,m}u^m,
 \qquad
 T_n(q)=a_n\E X_n^q=[z^n]G_q(z),
 \qquad
 G_q(z)=\prod_{k\geq1}g_q\!\left(\frac{z^k}{k^{2-q}}\right).
\end{equation}
For each coefficient, the last identity involves only finitely many factors
and terms. It therefore remains a valid formal identity when the series on
its right has radius of convergence zero.

\subsection{The moments which pass to the limit}

\begin{theorem}\label{thm:fractional}
For every fixed real number \(q<1\),
\begin{equation}\label{eq:fractional-moments}
 \E\left(\frac{X_n}{n}\right)^q
 \longrightarrow \E W^q=\frac{G_q(1)}{C}<\infty.
\end{equation}
For \(q\geq1\), \(\E W^q=\infty\).
In particular, the finite positive moments of the limiting variable are
exactly those of orders less than one.
\end{theorem}

\begin{proof}
We first treat \(0<q<1\), where convergence in distribution alone is
insufficient. Set \(\alpha=2-q\), so that \(1<\alpha<2\).
The ratio
\[
 \frac{v_{q,m+1}}{v_{q,m}}
 =\frac{(2m+1)^{q-1}}{2m+2}
 \longrightarrow0
\]
shows that \(g_q\) is entire. Near the origin,
\[
 g_q(u)=1+\frac{u}{2}+O(u^2),\qquad
 \log g_q(u)=\frac{u}{2}+\sum_{r\geq2}e_{q,r}u^r.
\]
In particular, \(G_q(1)\) is finite and positive: the logarithm of its
\(k\)-th factor is \(1/(2k^\alpha)+O(k^{-2\alpha})\).

The difficulty is that \(G_q\) can have singularities at infinitely many
roots of unity. We separate its principal singularity from a factor with
sufficiently summable coefficients, following the same division of labor
between singularity analysis and boundary regularity as in the hybrid
method \cite{FFGPP}. Choose
\[
 \alpha<\beta<2\alpha-1.
\]
Let \(\rho>0\) be small enough that the Taylor series for \(\log g_q\)
is analytic on a neighborhood of \(|u|\leq\rho\). Cauchy's estimate gives
\(|e_{q,r}|\leq M\rho^{-r}\) for a constant \(M\). Choose an integer
\(K\) with \(K^{-\alpha}<\rho\). Then
\begin{equation}\label{eq:fractional-hybrid}
 G_q(z)=P_\alpha(z)H_q(z),\qquad
 P_\alpha(z)=\exp\!\left(\frac12\Li_\alpha(z)\right),
\end{equation}
where
\begin{align*}
 H_q(z)&=
 \left\{\prod_{k\leq K}g_q\!\left(\frac{z^k}{k^\alpha}\right)\right\}
 \exp\!\left(-\frac12\sum_{k\leq K}\frac{z^k}{k^\alpha}\right)
 \exp R_K(z),\\
 R_K(z)&=\sum_{k>K}\sum_{r\geq2}
          \frac{e_{q,r}z^{kr}}{k^{\alpha r}}.
\end{align*}
No assertion about the zeros of the finitely many small-index factors is
needed; those factors have been left intact.

We claim that, if \(H_q(z)=\sum_{j\geq0}h_jz^j\), then
\begin{equation}\label{eq:weighted-summability}
 \sum_{j\geq0}(1+j)^\beta |h_j|<\infty.
\end{equation}
Here are the details. Since \(\beta-2\alpha<-1\),
\begin{align*}
 \sum_{k>K}\sum_{r\geq2}(1+kr)^\beta
       \frac{|e_{q,r}|}{k^{\alpha r}}
 &\leq C_\beta M\sum_{r\geq2}r^\beta\rho^{-r}
                  \sum_{k>K}k^{\beta-\alpha r}\\
 &\leq C_\beta M K^{\beta+1}
       \sum_{r\geq2}
       \frac{r^\beta}{\alpha r-\beta-1}
       \left(\frac{K^{-\alpha}}{\rho}\right)^r
 <\infty.
\end{align*}
Thus the coefficients of \(R_K\) have a finite weighted absolute sum.
For a series \(U(z)=\sum u_jz^j\), put
\(\|U\|_\beta=\sum(1+j)^\beta|u_j|\).
The inequality
\((1+i+j)^\beta\leq(1+i)^\beta(1+j)^\beta\)
shows that \(\|UV\|_\beta\leq\|U\|_\beta\|V\|_\beta\).
Consequently \(\|\exp R_K\|_\beta\leq\exp\|R_K\|_\beta\).
The other factor in \(H_q\) is entire, so Cauchy's estimate on any circle
of radius greater than one gives it a finite \(\beta\)-norm. This proves
\eqref{eq:weighted-summability}.

It remains to extract the coefficients of \(P_\alpha\). In a slit
neighborhood of one, with \(t=-\log z\), the polylogarithm expansion is
\[
 \Li_\alpha(z)
 =\zeta(\alpha)+\Gamma(1-\alpha)t^{\alpha-1}+O(t).
\]
Choose \(\eta>0\) smaller than both \(\alpha-1\) and \(2-\alpha\).
Since \(t=(1-z)+O((1-z)^2)\), exponentiation gives
\[
 P_\alpha(z)=e^{\zeta(\alpha)/2}
 \left\{1+\frac{\Gamma(1-\alpha)}2(1-z)^{\alpha-1}
       +O\!\left(|1-z|^{\alpha-1+\eta}\right)\right\}.
\]
The function \(P_\alpha\) is analytic in a standard indented domain at
one and across the rest of the unit circle. The transfer theorem
\cite{FO} therefore gives, for \(p_n=[z^n]P_\alpha(z)\),
\begin{equation}\label{eq:fractional-principal-coefficient}
 p_n\sim\frac{e^{\zeta(\alpha)/2}}{2n^\alpha},
 \qquad p_n=O((n+1)^{-\alpha}).
\end{equation}
The gamma factors cancel in the leading coefficient.

For clarity, the passage through the second factor in
\eqref{eq:fractional-hybrid} is justified directly. In
\(\sum_{j=0}^nh_jp_{n-j}\), split at \(j=n/2\).
On the first part, \(p_{n-j}/p_n\) is uniformly bounded and tends to
one for each fixed \(j\). Dominated convergence and
\eqref{eq:weighted-summability} show that this part, divided by
\(p_n\), tends to \(H_q(1)\). On the second part, the bound
\(p_{n-j}=O(1)\) gives
\[
 \frac1{p_n}\sum_{j>n/2}|h_jp_{n-j}|
 \leq O(n^\alpha)\sum_{j>n/2}|h_j|
 =o(n^{\alpha-\beta})=o(1).
\]
It follows that
\begin{equation}\label{eq:fractional-numerator}
 [z^n]G_q(z)\sim\frac{G_q(1)}{2n^\alpha}.
\end{equation}
Dividing by \(a_n\sim C/(2n^2)\) proves
\eqref{eq:fractional-moments} for \(0<q<1\).

The distribution of $\Omega$ in \eqref{eq:residual-mass} gives, by summation of nonnegative terms,
\[
 \E W^q
 =\frac1C\sum_{\omega}w(\omega)
      \prod_{k\geq1}\bigl(k^{m_k}r_{m_k}\bigr)^q
 =\frac{G_q(1)}C.
\]
For \(q<0\), the lower bound \(X_n\geq n-1\) for \(n\geq3\)
puts \(X_n/n\) in \([2/3,\infty)\), where the function
\(x\mapsto x^q\) is bounded and continuous. The distributional limit
therefore proves convergence of these negative moments. Their finiteness
also follows from \(W\geq1\). The case \(q=0\) is immediate.
Finally, for each \(k\geq1\), consider the event that \(M_k=1\) and
\(M_j=0\) for every \(j\ne k\). These events are disjoint and contribute
\[
 \E W^q\geq\frac1{2C}\sum_{k\geq1}k^{q-2}=\infty
 \qquad(q\geq1).
\]
\end{proof}

\subsection{The intermediate moments}

The interval $1<q<2$ lies between two different sources of large
moments. The limiting product $W$ already has an infinite $q$-th
moment, but the identity permutation has not yet taken over the moment
sum. In this interval an increasing collection of cycle lengths
contributes. Its generating function has an exponential singularity,
and a saddle-point estimate describes the resulting scale.

Fix $1<q<2$ and put $\alpha=2-q$ and $\beta=q-1=1-\alpha$. Throughout
this subsection these parameters are fixed. Define
\[
 \Phi_q(t)=\log G_q(e^{-t})\qquad(t>0).
\]
The logarithm here is real. The ratio of successive coefficients of
$g_q$ tends to zero, as in the proof of Theorem~\ref{thm:fractional};
thus $g_q$ is entire. Since $\log g_q(u)=u/2+O(u^2)$ at the origin,
the product defining $G_q$ converges for $|z|<1$.

\begin{theorem}\label{thm:intermediate}
For every positive integer $n$ there is a unique $t_n>0$ satisfying
\begin{equation}\label{eq:intermediate-saddle}
 -\Phi_q'(t_n)=n.
\end{equation}
As $n\to\infty$,
\begin{equation}\label{eq:intermediate-equivalent}
 \E X_n^q\sim\frac{2n^2}{C}
 \frac{\exp\{\Phi_q(t_n)+nt_n\}}
      {\sqrt{2\pi\Phi_q''(t_n)}}.
\end{equation}
Moreover,
\begin{equation}\label{eq:intermediate-scale}
 t_n\sim\left(\frac{\Gamma(q)}{2n}\right)^{1/q},\qquad
 \log\E X_n^q\sim
 \frac{q}{q-1}\left(\frac{\Gamma(q)}2\right)^{1/q}
 n^{(q-1)/q}.
\end{equation}
\end{theorem}

\begin{proof}
We first establish the estimates needed near the positive real axis.
Choose $K$ large enough that the series for $\log g_q(u)$ converges
absolutely on a disk containing $|u|\leq K^{-\alpha}$. Keep the
finitely many factors with $k\leq K$ intact. For the remaining
factors, logarithmic differentiation of the Taylor series gives,
for each fixed integer $j\geq0$,
\[
 \left(u\frac{\dd}{\dd u}\right)^j\log g_q(u)
 =\frac u2+O(u^2)
\]
uniformly on a sufficiently small disk. Consequently, for
$j=0,1,2,3$,
\begin{equation}\label{eq:intermediate-derivatives}
 (-1)^j\Phi_q^{(j)}(t)
 \sim\frac{\Gamma(\beta+j)}2t^{-\beta-j}.
\end{equation}
Indeed, the leading sum is
$\frac12\sum_{k\geq1}k^{j-\alpha}e^{-kt}$, whose equivalent follows
by a Riemann sum. The error is bounded by a constant times
$1+\sum_{k>K}k^{j-2\alpha}e^{-2kt}$. For $j\geq1$ this is
$O(t^{2\alpha-j-1})$, smaller than the main term by a factor
$t^\alpha$. For $j=0$ it is $O(t^{2\alpha-1})$ if
$2\alpha<1$, $O(\log(1/t))$ if $2\alpha=1$, and $O(1)$ otherwise.
Each is $o(t^{-\beta})$.

For $t>0$, assign independent nonnegative integer variables $Y_k$
the probabilities
\[
 \P_t(Y_k=m)=
 \frac{v_{q,m}(e^{-kt}/k^\alpha)^m}
      {g_q(e^{-kt}/k^\alpha)}.
\]
Their weighted sum $Y=\sum kY_k$ is almost surely finite, since its
expectation is finite. Differentiation gives
$\E_tY=-\Phi_q'(t)$ and $\Var_tY=\Phi_q''(t)>0$.
The expectation is strictly decreasing, tends to zero as
$t\to\infty$, and tends to infinity as $t\downarrow0$ by
\eqref{eq:intermediate-derivatives}. This proves existence and
uniqueness in \eqref{eq:intermediate-saddle}, as well as the first
equivalent in \eqref{eq:intermediate-scale}.

We next control the complementary arcs. Set $u_k=e^{-kt}/k^\alpha$
and let $p_{k,m}=\P_t(Y_k=m)$. Since $0<u_k\leq1$,
\[
 p_{k,0}p_{k,1}
 =\frac{u_k}{2g_q(u_k)^2}
 \geq\frac{u_k}{2g_q(1)^2}.
\]
Keeping only the terms with indices $0,1$ in the squared modulus
of the characteristic function therefore gives
\begin{equation}\label{eq:intermediate-modulus}
 \left|\frac{G_q(e^{-t+i\theta})}{G_q(e^{-t})}\right|
 \leq\exp\left\{-c_q\sum_{k\geq1}
  k^{-\alpha}e^{-kt}(1-\cos k\theta)\right\},
 \qquad |\theta|\leq\pi,
\end{equation}
for a positive constant $c_q$.

For completeness, the sum in \eqref{eq:intermediate-modulus}
satisfies the uniform lower bound
\begin{equation}\label{eq:intermediate-arc-bound}
 \sum_{k\geq1}k^{-\alpha}e^{-kt}(1-\cos k\theta)
 \geq c\,t^{-\beta}\min\{1,(\theta/t)^2\}
\end{equation}
for all sufficiently small $t>0$ and $|\theta|\leq\pi$.
To see this, put $N=\lfloor1/t\rfloor$ and restrict to $k\leq N$,
where the weights are at least $e^{-1}N^{-\alpha}$. If
$|\theta|\leq1/N$, the inequality $1-\cos x\geq c_0x^2$ for
$|x|\leq1$ gives
$\sum_{k\leq N}(1-\cos k\theta)\geq cN^3\theta^2$.
For $1/N\leq|\theta|\leq4\pi/N$, use instead
$k\leq\lfloor N/(4\pi)\rfloor$ to obtain a lower bound $cN$.
Finally, if $4\pi/N\leq|\theta|\leq\pi$, the geometric sum gives
\[
 \left|\sum_{k=1}^Ne^{ik\theta}\right|
 \leq\frac1{\sin(|\theta|/2)}
 \leq\frac\pi{|\theta|}\leq\frac N4,
\]
so the sum of $1-\cos k\theta$ is at least $3N/4$.
These three estimates prove \eqref{eq:intermediate-arc-bound}.

Now take $t=t_n$ and write $B=\Phi_q''(t)$. On a neighborhood of
the central arc the logarithm can be chosen continuously from its
value at $\theta=0$. No assertion about zeros elsewhere is needed:
the finitely many small-index factors are nonzero near their
positive real arguments, and the remaining logarithms are
represented by their convergent Taylor series. Those series also
give the uniform bound
\[
 |\Phi_q'''(t-iv)|=O(t^{-\beta-3})
 \qquad(|v|\leq t).
\]
Fix $0<\varepsilon<\beta/6$ and set
$\theta_0=B^{-1/2}t^{-\varepsilon}$. Then
$\theta_0/t\to0$, and Taylor's formula gives, uniformly on
$|\theta|\leq\theta_0$,
\[
 \Phi_q(t-i\theta)-\Phi_q(t)-in\theta
 =-\frac12B\theta^2+O(t^{-\beta-3}|\theta|^3)
 =-\frac12B\theta^2+o(1).
\]
Here the last error is $O(t^{\beta/2-3\varepsilon})$.
The central integral is consequently asymptotic to
$\sqrt{2\pi/B}$. By \eqref{eq:intermediate-modulus} and
\eqref{eq:intermediate-arc-bound}, the part with
$\theta_0\leq|\theta|\leq t$ is bounded by a Gaussian tail and is
$o(B^{-1/2})$; the remaining part is
$O(\exp(-ct^{-\beta}))=o(B^{-1/2})$.
Cauchy's formula therefore yields
\[
 T_n(q)\sim
 \frac{\exp\{\Phi_q(t_n)+nt_n\}}
      {\sqrt{2\pi\Phi_q''(t_n)}}.
\]
Division by $a_n\sim C/(2n^2)$ proves
\eqref{eq:intermediate-equivalent}. Finally,
\eqref{eq:intermediate-derivatives} and
$t_n\sim(\Gamma(q)/(2n))^{1/q}$ give the logarithmic equivalent
in \eqref{eq:intermediate-scale}; the logarithm of the prefactor
is only $O(\log n)$.
\end{proof}

The saddle in Theorem~\ref{thm:intermediate} retains all lower-order
singular terms automatically. A finite version makes their role
more explicit. Write
\[
 \log g_q(u)=\sum_{r\geq1}e_{q,r}u^r,
 \qquad e_{q,1}=\frac12,\quad
 e_{q,2}=\frac{3^q-3}{24},
 \qquad J=\lfloor1/\alpha\rfloor.
\]
For $1\leq r\leq J$ with $\alpha r<1$, put
\[
 p_r=1-\alpha r,\qquad
 A_r=e_{q,r}\Gamma(p_r)r^{-p_r}.
\]
Set $b=e_{q,J}$ if $\alpha J=1$, and $b=0$ otherwise. Define
\begin{align}
 \Psi_q(t)&=\sum_{\substack{1\leq r\leq J\\\alpha r<1}}
                 A_rt^{-p_r}+b\log(1/t),\label{eq:finite-saddle-function}\\
 D_q&=\sum_{k\geq1}\left\{
   \log g_q(k^{-\alpha})-
   \sum_{r=1}^Je_{q,r}k^{-\alpha r}\right\}
   +\sum_{\substack{1\leq r\leq J\\\alpha r<1}}
           e_{q,r}\zeta(\alpha r)-b\log J.
 \label{eq:intermediate-constant}
\end{align}
The first sum in \eqref{eq:intermediate-constant} converges
absolutely, since its summand is $O(k^{-\alpha(J+1)})$.
The zeta values below one mean analytic continuation, not
convergent Dirichlet series.

\begin{corollary}\label{cor:finite-intermediate-saddle}
For all sufficiently large $n$, the equation
\[
 -\Psi_q'(\widehat t_n)
 =\sum_{\alpha r<1}p_rA_r\widehat t_n^{-p_r-1}
       +b/\widehat t_n=n
\]
has a unique solution in a fixed sufficiently small interval
$(0,t_0)$. With that solution,
\begin{equation}\label{eq:finite-intermediate-equivalent}
 \E X_n^q\sim
 \frac{2e^{D_q}n^2}{C}
 \frac{\exp\{n\widehat t_n+\Psi_q(\widehat t_n)\}}
      {\sqrt{2\pi\Psi_q''(\widehat t_n)}}.
\end{equation}
\end{corollary}
\begin{proof}
Put $\delta=\alpha(J+1)-1$, so $0<\delta<1$.
Subtracting the first $J$ terms of $\log g_q$ gives
\[
 \Phi_q(t)=\sum_{r=1}^Je_{q,r}\Li_{\alpha r}(e^{-rt})
       +\mathcal R_q(t).
\]
The Taylor remainder and its logarithmic derivatives have order
$u^{J+1}$ at zero. Summation consequently gives
\[
 \mathcal R_q(t)=\mathcal R_q(0)+O(t^\delta),\qquad
 \mathcal R_q^{(j)}(t)=O(t^{\delta-j})\quad(j=1,2).
\]
For the first estimate, bound the summands of
$\mathcal R_q(t)-\mathcal R_q(0)$ by
$Ck^{-1-\delta}\min(kt,1)$ and split at $k=1/t$.
The derivative estimates follow from
$\sum k^{j-1-\delta}e^{-(J+1)kt}=O(t^{\delta-j})$.
The polylogarithm expansions at the positive real axis now give
\begin{equation}\label{eq:intermediate-finite-expansion}
 \Phi_q(t)=\Psi_q(t)+D_q+O(t^\delta),\qquad
 \Phi_q^{(j)}(t)=\Psi_q^{(j)}(t)+O(t^{\delta-j})
 \quad(j=1,2).
\end{equation}
In particular, a term with $\alpha r=1$ contributes
$e_{q,r}\log(1/t)-e_{q,r}\log r$, which explains the last term
of \eqref{eq:intermediate-constant}.

The leading term of $\Psi_q$ is
$\Gamma(\beta)t^{-\beta}/2$, so $-\Psi_q'$ is strictly
decreasing from infinity on a sufficiently small interval. This
proves the stated existence and uniqueness. Both saddle points
are asymptotic to $(\Gamma(q)/(2n))^{1/q}$, and the mean value
theorem applied to their defining equations gives
\[
 t_n-\widehat t_n=O(\widehat t_n^{\delta+\beta+1}).
\]
The change in the saddle exponent is therefore
$O(\widehat t_n^{2\delta+\beta})=o(1)$, while the change in
the variance is relative $o(1)$. Substitution in
Theorem~\ref{thm:intermediate}, using
\eqref{eq:intermediate-finite-expansion}, proves the result.
\end{proof}

For example, if $1<q<3/2$, only the term $r=1$ is singular.
Writing $\lambda_q=\Gamma(q)/2$ and
$K_q=q\lambda_q^{1/q}/(q-1)$, the corollary simplifies to
\begin{equation}\label{eq:intermediate-first-range}
 \E X_n^q\sim
 \frac{2e^{D_q}\lambda_q^{1/(2q)}}{C\sqrt{2\pi q}}
 n^{3/2-1/(2q)}\exp\{K_qn^{(q-1)/q}\}.
\end{equation}
At $q=3/2$, the $r=2$ term becomes logarithmic; for larger $q$
it supplies an additional stretched exponential. Further terms
can enter when $q=2-1/r$. Thus the leading logarithmic scale is
simple throughout $1<q<2$, but a relative asymptotic requires
the additional terms recorded in \eqref{eq:finite-saddle-function}.
These estimates are for fixed $q$; no uniformity as $q$ tends to
$1$ or $2$ is asserted.

\subsection{Identity domination above order two}

The next result requires a bound on the sum of all nonidentity
contributions. Comparing only the largest summands would not suffice,
since the number of sequences of cycle counts grows with \(n\).

\begin{theorem}\label{thm:higher}
For every fixed real \(q>2\),
\begin{equation}\label{eq:high-moment-asymptotic}
 \E X_n^q\sim
 \frac{2^{q/2}}{C\sqrt\pi}\,
 n^{3/2}\left(\frac{2n}{e}\right)^{(q-2)n}.
\end{equation}
More precisely, putting \(p=q-2\), the moment numerator satisfies
\begin{equation}\label{eq:identity-domination}
 \frac{T_n(q)}{v_{q,n}}
 =1+2^{-p-1}n^{-2p}
       +O\!\left(n^{-2p-\min(p,1)}\right).
\end{equation}
Thus the identity permutation supplies an asymptotic proportion one of
the \(q\)-th moment numerator.
\end{theorem}

\begin{proof}
Recall \(d_m=\binom{2m}{m}/4^m\), so that
\(v_{q,m}=d_mr_m^p\). Separate a sequence of cycle counts according to
\(j=n-m_1\), half the number of letters in its cycles of length at least two. Define
\[
 U_j(q)=\sum_{\substack{\sum_{k\geq2}km_k=j\\m_k\geq0}}
          \prod_{k\geq2}k^{pm_k}v_{q,m_k}.
\]
Then
\begin{equation}\label{eq:nonfixed-size-split}
 T_n(q)=\sum_{j=0}^n v_{q,n-j}U_j(q),\qquad
 U_0(q)=1,\quad U_1(q)=0,\quad U_2(q)=2^{p-1}.
\end{equation}

Consider a sequence \((m_k)_{k\geq2}\) contributing to \(U_j\), so that the nonfixed cycles contain \(2j\) letters, and set
\(M=\sum_{k\geq2}m_k\). Since \(M\leq j/2\) and
\(r_m\leq(2m)^m\), its factorization count is at most
\[
 \prod_{k\geq2}k^{m_k}r_{m_k}
 \leq\prod_{k\geq2}(2km_k)^{m_k}
 \leq(2j)^M\leq(2j)^{j/2}.
\]
On the other hand, the sum of its second-moment weights is
\[
 [z^j]\prod_{k\geq2}(1-z^k)^{-1/2}
 \leq [z^j]\prod_{k\geq1}(1-z^k)^{-1}
 \leq 2^j.
\]
The first inequality is coefficientwise, since
\(\binom{2m}{m}/4^m\leq1\); the second bounds partitions of
\(j\) by compositions. Multiplying a second-moment weight by the
\(p\)-th power of the factorization count proves
\begin{equation}\label{eq:residual-moment-bound}
 U_j(q)\leq2^j(2j)^{pj/2}\qquad(j\geq2).
\end{equation}

Suppose first that \(2\leq j\leq n/2\). Each of the \(j\) odd
factors omitted from \(r_n\) to obtain \(r_{n-j}\) is at least
\(n\), whence \(r_{n-j}/r_n\leq n^{-j}\).
Also \(d_{n-j}/d_n\) is bounded uniformly in this range, by the
central binomial coefficient estimate. Thus
\[
 \frac{v_{q,n-j}U_j(q)}{v_{q,n}}
 \leq C_p T_{j,n},\qquad
 T_{j,n}=2^j(2j)^{pj/2}n^{-pj}.
\]
These bounds are summable uniformly in \(n\), because
\[
 \frac{T_{j+1,n}}{T_{j,n}}
 =2n^{-p}[2(j+1)]^{p/2}
       \left(1+\frac1j\right)^{pj/2}
 \leq2[2e(j+1)]^{p/2}n^{-p}
 \leq C_p n^{-p/2}.
\]
For sufficiently large \(n\), the last expression is less than
\(1/2\). Consequently the sum beginning at \(j=2\) is
\(O(n^{-2p})\), and the sum beginning at \(j=3\) is
\(O(n^{-3p})\).

For \(j>n/2\), use instead
\(r_n\geq r_{n-j}r_j\) and
\(r_j\geq j!\geq(j/e)^j\). The central binomial coefficient
estimate now gives \(d_{n-j}/d_n=O(\sqrt n)\), uniformly even at
\(j=n\). It follows from \eqref{eq:residual-moment-bound} that
\[
 \frac{v_{q,n-j}U_j(q)}{v_{q,n}}
 \leq C_p\sqrt n\left\{2(e\sqrt2)^p j^{-p/2}\right\}^{j}.
\]
For some constant \(c_p>0\), the sum over \(n/2<j\leq n\)
is at most \(\exp(-c_p n\log n)\) for all sufficiently large \(n\). This proves that the nonidentity
contribution is \(O(n^{-2p})\) times the identity contribution,
and also justifies isolating the \(j=2\) term.

That term can be evaluated without approximation:
\[
 \frac{v_{q,n-2}U_2(q)}{v_{q,n}}
 =2^{p-1}\frac{(2n)(2n-2)}
                {\bigl((2n-1)(2n-3)\bigr)^{p+1}}
 =2^{-p-1}n^{-2p}\bigl(1+O(n^{-1})\bigr).
\]
Together with the \(O(n^{-3p})\) bound for the remaining terms,
this proves \eqref{eq:identity-domination}.
Finally, Stirling's formula gives
\[
 r_n=\sqrt2\left(\frac{2n}{e}\right)^n(1+O(n^{-1})),
 \qquad
 d_n=\frac1{\sqrt{\pi n}}(1+O(n^{-1})).
\]
Hence
\[
 T_n(q)\sim v_{q,n}
 \sim\frac{2^{(q-2)/2}}{\sqrt{\pi n}}
       \left(\frac{2n}{e}\right)^{(q-2)n}.
\]
Division by \(a_n\sim C/(2n^2)\) yields
\eqref{eq:high-moment-asymptotic}.
\end{proof}

The change at order one is a failure of uniform integrability on the
typical scale \(n\). At order two, a partition generating function
appears. Above order two, the factorial growth of the contribution from
fixed points outweighs every other sequence of cycle counts, including their total
number. Theorem~\ref{thm:intermediate} fills the intermediate range \(1<q<2\). All these statements concern fixed moment orders; no uniform assertion across the transition points is intended.

\section{Concluding remarks}
The three refinements suggested by the qualitative limit law can now be made precise. Theorem~\ref{thm:sharp-median} gives the sharp median threshold: the unique median is $n$ for every $n\geq4$, and the only exception to the formula among positive half-degrees is $n=3$. Theorem~\ref{thm:intermediate} supplies a full saddle-point equivalent for every fixed moment order between one and two. Theorem~\ref{thm:residual-rate} gives the sharp total variation expansion for the cycle counts left after removal of the two long cycles. Proposition~\ref{prop:no-giant-and-count} then uses these counts to approximate the factorization count itself.

For every $n$, the sequences of remaining cycle counts belong to the fixed countable set $\mathcal C$. Including the symbol $\ast$ allows us to compare both distributions on $\mathcal C\cup\{\ast\}$, with probability zero assigned to any value that cannot occur. Total variation makes this comparison by summing absolute differences of probabilities. For the intermediate moments, the first-order logarithm has a simple expression, while a relative equivalent must retain every contributing singular term. The finite saddle in Corollary~\ref{cor:finite-intermediate-saddle} does this without requiring a separate proof in each subinterval of moment orders.

The estimates here keep the moment order fixed. Uniform asymptotics when it approaches one or two with $n$ remain to be developed, as do explicit finite-degree error bounds for the saddle-point approximation. These are further refinements of the laws established above.

\section{Acknowledgments and Disclosure of AI Use}

This paper ties up all of the loose ends from the author's doctoral thesis \cite{BurnetteThesis}. I would like to thank my advisor, Dr. Eric Schmutz, for his continued support, both during my graduate studies and for the near decade after.

The ideas articulated in this paper are my own. However, due to laziness and a shortage of free time, I fed my fragmented notes to GPT-6.2 Astra (OpenAI) and asked it to synthesize everything into a completed manuscript. I take responsibility for the mathematical arguments, computations, and references in this paper.

\clearpage
\appendix
\section{An overview of singularity analysis}\label{app:singularity}
The coefficient argument in Section~\ref{sec:enumeration} merges two
different ways of controlling an error term. Darboux's method uses
smoothness on the unit circle, while singularity analysis uses analytic
continuation and the local size of a function near its singularities.
We collect the relevant principles here, with the normalization needed
for the acme permutation generating function. This also explains why the two
singularities of the full-degree series $E(z)=F(z^2)$ yield the same
asymptotic data as the half-degree series $F(z)$.

\subsection{Darboux's method}
In this subsection, we collect basic facts about functions of moderate
growth near the unit circle. Write $\mathbb D=\{z:|z|<1\}$ for the open unit disc.

\begin{definition}\label{def:global-order}
An analytic function $f$ in $\mathbb D$ has \emph{global order}
$\eta\leq0$ if $f(z)=O((1-|z|)^\eta)$ throughout $\mathbb D$.
\end{definition}
We are particularly interested in functions whose local expansions near dominant singularities
combine powers and logarithms. Put $\mathcal L(z)=\log(1/(1-z))$.

\begin{definition}\label{def:log-power}
A \emph{log-power function at $1$} is a finite sum of the form
\[
 \sigma(z)=\sum_{j=1}^v p_j(\mathcal L(z))(1-z)^{\alpha_j},
\]
where the $p_j$ are polynomials and $\alpha_1<\cdots<\alpha_v$ are real.
A log-power function at a finite set $Z$ on the unit circle $\partial\mathbb{D}$ is a sum
$\Sigma(z)=\sum_{\xi\in Z}\sigma_\xi(z/\xi)$ of such functions.
All branches agree with their power series in $\mathbb D$.
\end{definition}
The coefficients of these functions admit full expansions in decreasing
powers of the coefficient index; see Jungen~\cite{Jungen} and
\cite[Chapter~VI]{FS}. We record the terms and differentiation rule
needed here.

\begin{proposition}\label{prop:appendix-logpowers}
For fixed $\alpha\notin\{0,1,2,\ldots\}$,
\begin{equation}\label{eq:appendix-gamma-expansion}
 [z^d](1-z)^\alpha
 =\frac{d^{-\alpha-1}}{\Gamma(-\alpha)}
  \left(1+\frac{\alpha(\alpha+1)}{2d}+O(d^{-2})\right).
\end{equation}
For every nonnegative integer $q$, the exact differentiation rule is
\[
 [z^d](1-z)^\alpha\mathcal L(z)^q
 =(-1)^q\frac{\partial^q}{\partial\alpha^q}
   \left[\frac{\Gamma(d-\alpha)}{\Gamma(-\alpha)\Gamma(d+1)}\right].
\]
The quotient is interpreted by analytic continuation when necessary.
Consequently, for fixed $q$ and $\alpha\notin\{0,1,2,\ldots\}$,
\[
 [z^d](1-z)^\alpha\mathcal L(z)^q
 \sim\frac{d^{-\alpha-1}(\log d)^q}{\Gamma(-\alpha)}.
\]
For integers $r\geq0$ and $q\geq1$, the corresponding equivalent is
\begin{equation}\label{eq:appendix-integral-logpower}
 [z^d](1-z)^r\mathcal L(z)^q
 \sim(-1)^r q\,r!\,d^{-r-1}(\log d)^{q-1}.
\end{equation}
\end{proposition}

For example, for $d>1$,
\[
 [z^d]\mathcal L(z)=\frac1d,\qquad
 [z^d](1-z)\mathcal L(z)=-\frac1{d(d-1)}.
\]
Thus a logarithm multiplied by a vanishing factor can still determine
the first nonzero coefficient equivalent. At any other point $\xi$,
the substitution rule
$[z^d]\psi(z/\xi)=\xi^{-d}[z^d]\psi(z)$ accounts for the oscillation
in the coefficients.

\begin{definition}\label{def:smooth-expansion}
A function $f$, analytic in $\mathbb D$, is \emph{$C^s$-smooth} on the
closed disc $\overline{\mathbb{D}}$ if each derivative of order at most the nonnegative integer
$s$ extends continuously to $\overline{\mathbb D}$. It admits a
\emph{log-power expansion of class $C^s$} if $f-\Sigma$ is $C^s$-smooth
for some log-power function $\Sigma$ at finitely many points on $\partial\mathbb{D}$.
\end{definition}

\begin{proposition}[Based on Olver~\cite{Olver}, Theorem~8.9.1]\label{prop:appendix-darboux}
If $f-\Sigma$ is $C^s$-smooth on $\overline{\mathbb{D}}$, then $[z^d]f(z)=[z^d]\Sigma(z)+o(d^{-s}).$
\end{proposition}

The smoothness belongs to the remainder. Our $E(z)$ and $F(z)$ is Section 2 have global
order 0, since their coefficients are nonnegative and their values
at $z=1$ equal $C$. Nevertheless, the leading logarithmic term in
\eqref{eq:F-local} shows that they are not $C^1$ at $z=1$.

\subsection{The Flajolet--Odlyzko transfer theorem}
For $R_0>1$ and $0<\phi<\pi/2$, define
\[
 \Delta(R_0,\phi)=
 \{z:|z|<R_0,\ z\ne1,\ |\arg(z-1)|>\phi\},
\]
where the argument is taken in $(-\pi,\pi]$. This domain contains
$\mathbb D$ and excludes a wedge pointing outward from $1$.
For a finite set $Z = \{\xi_1, \ldots, \xi_m\}$ on $\partial\mathbb{D}$, an intersection of rotated
domains $ \mathfrak{D} = \bigcap_{j=1}^{m} \xi_j\Delta(R_0,\phi)$, $\xi\in Z$, supplies the corresponding
indentations at all points of $Z$.

\begin{definition}\label{def:transfer-expansion}
A function admits a \emph{log-power expansion of type $O_t$} at $Z$
if it continues analytically to such an indented domain $\mathfrak{D}$ and, for every
$\xi\in Z$, has an expansion
\[
 f(z)=\sigma_\xi(z/\xi)+O(|z-\xi|^t)
\]
as $z\to\xi$ within $\xi\Delta(R_0,\phi),$ where $\sigma_\xi$ is a log-power function
at $1$.
\end{definition}

\begin{proposition}[Flajolet--Odlyzko \cite{FO}]
\label{prop:appendix-transfer}
For a function with a log-power expansion of type $O_t$ at a finite set
$Z$ on the unit circle,
\begin{equation}\label{eq:appendix-transfer}
 [z^d]f(z)=\sum_{\xi\in Z}[z^d]\sigma_\xi(z/\xi)+O(d^{-t-1}).
\end{equation}
\end{proposition}
This is a transfer theorem because it translates a local bound near
each singularity into a coefficient bound. Its hypothesis includes
analytic continuation, which is absent from a merely radial expansion.
Analytic polynomial terms in these local expansions have no effect
on sufficiently high coefficients. They must still be retained when
expansions are multiplied: an analytic linear term can multiply a
logarithmic term and affect the next coefficient correction. This
is precisely how $h=R'(1)$ enters \eqref{eq:an-two-terms}.

\subsection{Hybridization of Darboux's method and the transfer theorem}
The hybrid method allows us to work with radial expansions even when
the entire function does not continue across the unit circle. The
continuation is required of only one factor; smoothness controls the
other factor.

\begin{definition}\label{def:radial-expansion}
A \emph{radial expansion of order $t$} at $\xi\in\partial\mathbb D$
is a log-power function $\sigma_\xi(z/\xi)$ satisfying
\[
 f((1-x)\xi)=\sigma_\xi(1-x)+O(x^t)\qquad(x\downarrow0).
\]
We discard terms individually of order $O(x^t)$ and denote the
remaining expansion by $\operatorname{asymp}(f,\xi,t)$.
\end{definition}

\begin{proposition}[Hybrid method, {\cite[Theorem~2]{FFGPP}}]
\label{prop:appendix-hybrid}
Suppose $f=PQ$ is analytic in $\mathbb D$ and has finite global order
$\eta\leq0$. Let $Z$ be a finite subset of $\partial\mathbb{D}$. Assume that
$Q$ is $C^s$-smooth for a nonnegative integer $s$, and that $P$ has a
log-power expansion of type $O_{t_0}$ at $Z$, with $t_0\geq0$. Put
\[
 u_0=\left\lfloor\frac{s+\lfloor\eta\rfloor}{2}\right\rfloor.
\]
If $t_0>u_0$, then $f$ has radial expansions of order $u_0$ at the
points of $Z$, and
\[
 [z^d]f(z)=
 [z^d]\sum_{\xi\in Z}\operatorname{asymp}(f,\xi,u_0)
 +o(d^{-u_0}).
\]
\end{proposition}
Thus the method reduces coefficient extraction to two tasks: exhibit
an appropriate factorization $PQ$, and compute the radial expansions
at its finitely many relevant boundary points. The global-order
parameter enters as $\lfloor\eta\rfloor$, not as $|\eta|$.
In the bounded case used in this paper, $\eta=0$. Taking $s=4$ gives
$u_0=2$, which is sufficient for a main term of order $d^{-2}$.
Taking only $s=2$ would give an error $o(d^{-1})$, which is insufficient.

Lemma~\ref{lem:hybrid} proves the bounded version required here directly.
Its proof uses Hermite interpolation to replace the smooth factor by a
polynomial matching its derivatives at the exceptional points. Darboux's
method controls the difference, and the transfer theorem applies to the
remaining finite collection of singular expansions.

\subsection{Polylogarithms and the two dominant singularities}
Polylogarithms make this separation explicit. For $\nu\in\mathbb C$,
\[
 \Li_\nu(z)=\sum_{j\geq1}\frac{z^j}{j^\nu}\qquad(|z|<1).
\]
They continue analytically to $\mathbb C\setminus[1,\infty)$.
For a positive integer $m$, the convergent local expansion is
\begin{equation}\label{eq:appendix-polylog-integer}
 \Li_m(e^{-\tau})
 =\frac{(-\tau)^{m-1}}{(m-1)!}
      (H_{m-1}-\log\tau)
  +\sum_{\substack{j\geq0\\j\ne m-1}}
        \zeta(m-j)\frac{(-\tau)^j}{j!},
\end{equation}
where $H_0=0$, $H_r=\sum_{j=1}^rj^{-1}$, and
\[
 \tau=-\log z=\sum_{j=1}^{\infty}\frac{(1-z)^j}{j}.
\]
For $\nu\notin\{1,2,3,\ldots\}$, the corresponding expansion is
\begin{equation}\label{eq:appendix-polylog-noninteger}
 \Li_\nu(e^{-\tau})
 =\Gamma(1-\nu)\tau^{\nu-1}
     +\sum_{j=0}^{\infty}\zeta(\nu-j)\frac{(-\tau)^j}{j!}.
\end{equation}
Here $|\tau|<2\pi$ and a consistent branch is chosen in the slit
neighborhood; see \cite[p.~387]{Lewin} and \cite[Chapter~VI]{FS}.
In particular, the first nonanalytic term of $\Li_m(z)$ for $m\geq2$
has order $(1-z)^{m-1}\log(1-z)$. Direct differentiation of the
power series shows that $\Li_m$ is $C^{m-2}$ on $\overline{\mathbb{D}}$,
while its next derivative has a logarithmic singularity at $1$.

A useful regularity observation must include a convergence hypothesis.
If $r(\ell)$ is polynomially bounded and $k\geq2$, then the centered
series
\[
 \sum_{\ell= k}^{\infty}r(\ell)
   \{\Li_\ell(z)-z-[\zeta(\ell)-1]\}
\]
is $C^{k-2}$-smooth. Indeed, after the $z$ term has been removed,
the differentiated coefficients start at index $2$, and their absolute
sum is $O(2^{-\ell})$ for each derivative of order at most $k-2$.
Polynomial growth in $\ell$ is therefore harmless. Without removal
of the $z$ term, polynomial boundedness alone would not suffice:
$\Li_\ell(z)-\zeta(\ell)$ tends to $z-1$ as $\ell\to\infty$.
For our exp-log series, the stronger geometric bound
\eqref{eq:bj-bound} supplies absolute convergence directly.

For the full-degree series, write, for each nonnegative integer $J,$
\[
 E(z)=\widetilde P_J(z)\widetilde Q_J(z),\qquad
 \widetilde P_J(z)=
 \exp\left(\sum_{j=1}^Jb_j\Li_{2j}(z^{2j})\right),\qquad
 \widetilde Q_J(z)=
 \exp\left(\sum_{j=J}^{\infty}b_j\Li_{2j}(z^{2j})\right).
\]
For $0\leq r\leq2J$ and $j>J$, direct differentiation gives
\[
 \sup_{|z|\leq1}
 \left|\frac{\dd^r}{\dd z^r}\Li_{2j}(z^{2j})\right|
 \leq(2j)^r\zeta(2j-r).
\]
The smallest zeta argument is $2$. Hence the differentiated tail
converges uniformly and $\widetilde Q_J$ is $C^{2J}$. In general
one cannot replace $2J$ by $2J+2$: the first retained term of the
tail already has a logarithmic singularity in derivative $2J+1$.
Taking $J=2$ gives the $o(d^{-2})$ remainder needed for the leading
enumeration term. Higher precision requires retaining more
polylogarithms, as in Section~\ref{sec:enumeration}.

Finally, substitute $z^2$ for $z$ in \eqref{eq:F-local}. As $z\to1$
radially, with $u=1-z$, this yields
\[
 \frac{E(z)}C
 =1+u\log u+(\log2-1-2h)u+O(u^2(\log u)^2).
\]
Since $E(-z)=E(z)$, the same expansion holds at $-1$ with $u=1+z$.
The hybrid argument just given permits coefficient extraction from
these local expansions. The leading singular terms therefore give
\[
 [z^d]E(z)
 =\frac{C(1+(-1)^d)}{d(d-1)}+o(d^{-2}).
\]
In odd degree the coefficient is exactly zero by symmetry. In
degree $d=2n$, the two contributions add and give
$[z^{2n}]E(z)\sim C/(2n^2)$, in agreement with
Theorem~\ref{thm:enumeration}. The half-degree variable packages
these two equal contributions into a single dominant singularity.

\end{document}